\documentclass[11pt]{article}
\usepackage{amsmath}
\numberwithin{equation}{section}
\usepackage{amsfonts,amssymb,amsmath,amsthm}
\usepackage{mathrsfs}
\newcommand{\beq}{\begin{equation}}
\newcommand{\enq}{\end{equation}}

\newtheorem{Theorem}{Theorem}[section]
\newtheorem{Lemma}[Theorem]{Lemma}

\newtheorem{Definition}[Theorem]{Definition}
\newtheorem{Remark}[subsection]{Remark}

\newcommand{\benu}{\begin{enumerate}}
\newcommand{\beqa}{\begin{eqnarray}}
\newcommand{\beqan}{\begin{eqnarray*}}
\newcommand{\eay}{\end{array}}
\newcommand{\edm}{\end{displaymath}}
\newcommand{\eenu}{\end{enumerate}}
\newcommand{\eeq}{\end{equation}}
\newcommand{\eeqa}{\end{eqnarray}}
\newcommand{\eeqan}{\end{eqnarray*}}

\newcommand{\br}{\begin{Remark}}
\newcommand{\er}{\end{Remark}}

\newcommand{\bqa}{\begin{eqnarray}}
\newcommand{\eqa}{\end{eqnarray}}
\newcommand{\bqw}{\begin{eqnarray*}}
\newcommand{\eqw}{\end{eqnarray*}}

\newcommand{\bea}{\begin{array}{cc}}
\newcommand{\ena}{\end{array}}

\allowdisplaybreaks[4]
\begin{document}
\begin{center}

{\large \bf  Upper-semicontinuity of uniform attractors for the non-autonomous viscoelastic Kirchhoff plate equation with memory }\\

 \vspace{0.20in}Yuming Qin$^{1,2\ast}$ $\ $ Hongli Wang$^{2}$\\
\end{center}
 $^{1}$ Insitute of Nonlinear Science, Donghua University,
Shanghai, 201620,  P. R. China, \\
$^{2}$ School of Mathematics and Statistics, Donghua University,
Shanghai, 201620, P. R. China.
 \vspace{3mm}

\begin{abstract}
	This paper delves into the long-time dynamics of a non-autonomous viscoelastic Kirchhoff plate equation with memory effects, described by 
	$$
		u_{t t}-\Delta u_{t t}+a_\epsilon(t) u_t+\alpha \Delta^2 u-\int_0^{\infty} \mu(s) \Delta^2 u(t-s) \mathrm{d} s-\Delta u_t+f(u)=g(x,t),
	$$
	in bounded domain $\Omega \subset \mathbb{R}^N$ with smooth boundary and nonlinear terms. Initially, the global existence of a weak solution that induces a continuous process is established. Subsequently, the existence of a uniform attractor is demonstrated in both subcritical and critical growth scenarios, utilizing operator techniques and an innovative analytical approach. Finally, the upper semicontinuity of the family of uniform attractors as the pert parameterurbation $\epsilon \to 0^+$ is proven through delicate energy estimates and a contradiction argument. Our results not only extend classical attractor theory to more general non-autonomous viscoelastic systems but also resolve open questions regarding the limiting behavior of attractors in the presence of both memory and critical nonlinearity.
 \vskip 3mm
 
 
 \smallskip
 \smallskip
 
\end{abstract}

\section{Introduction}
\setcounter{equation}{0}

\let\thefootnote\relax\footnote{*Corresponding author.}
\let\thefootnote\relax\footnote{E-mails: yuming@dhu.edu.cn (Y. Qin), whl@mail.dhu.edu.cn (H. Wang).}
\quad
In this paper, we are concerned with the following non-autonomous viscoelastic Kirchhoff plate equation with memory
\begin{equation}\label{eq1.1}
	u_{t t}-\Delta u_{t t}+a_\epsilon(t) u_t+\alpha \Delta^2 u-\int_0^{\infty} \mu(s) \Delta^2 u(t-s) \mathrm{d} s-\Delta u_t+f(u)=g(x,t),
\end{equation}
with simply supported boundary condition
\begin{equation}\label{eq1.2}
	u=\Delta u=0 \quad \text { on } \quad \partial \Omega \times \mathbb{R},
\end{equation}
and initial conditions
\begin{equation}\label{eq1.3}
	u(x, \tau)=u_0(x, \tau) \quad \text { and } \quad u_t(x, \tau)=\partial_t u_0(x, \tau), \quad(x, \tau) \in \Omega \times(-\infty, 0],
\end{equation}
where $\Omega$ is a bounded domain of $\mathbb{R}^N$ that has a smooth boundary, denoted by $\partial \Omega, u_0: \Omega \times(-\infty, 0] \rightarrow \mathbb{R}$ is defined as the given past history of $u$. The function $u_0: \Omega \times (-\infty, 0] \rightarrow \mathbb{R}$ prescribes the historical profile of the unknown $u$. The parameter $\alpha>0$  and $\mu$ represents a nonnegative kernel function associated with memory effects.
 
 The study of Kirchhoff-type equations has seen significant advancements, particularly in understanding their well-posedness and long-term dynamics. Kirchhoff's original work in 1983 introduced the canonical form of the Kirchhoff wave equation
 $$
 u_{t t}-M\left(\|\nabla u\|^2\right) \Delta u=h(x),
 $$ 
 which laid the foundation for subsequent research \cite{GK17}. Building on this, Chueshov \cite{IC18} explored the well-posedness and long-term behavior of an autonomous Kirchhoff wave equation with nonlinear damping
 $$
 u_{t t}-\phi\left(\|\nabla u\|^2\right) \Delta u-\sigma\left(\|\nabla u\|^2\right) \Delta u_t+f(u)=h(x),
 $$establishing the existence of a global attractor and its bounded Hausdorff and fractal dimensions under subcritical growth conditions. By looking at the trajectory of research in this field, it is evident that the majority of these studies, including those in \cite{IC23, MJ24, YQ22, ZY25}, has focused on clarifying the global well-posedness and intricate properties of global attractors for autonomous Kirchhoff equations. This focus reflects a sustained effort to deepen our understanding of the underlying mathematical structures and physical phenomena governed by these equations.
 
 In the realm of non-autonomous systems, Li et al. \cite{YL20} examined non-autonomous Kirchhoff wave equations with strong damping
 $$
 u_{t t}-\left(1+\epsilon\|\nabla u\|^2\right) \Delta u-\Delta u_t+f(u)=g(x, t),
 $$
focusing on the existence and upper-semicontinuity of uniform attractors. Yang et al. \cite{BY21} extended this research to time-dependent spaces, introducing complexities such as time-varying coefficients and evolving external forces $$
\varepsilon(t) \partial_{t t} u-\left(1+\delta\|\nabla u\|^2\right) \Delta u-\Delta \partial_t u+\lambda u=g(u)+h(t).
$$
Their work rigorously proved the existence and upper semi-continuity of pullback attractors.

 Sun, Cao and Duan \cite{CS6} investigated non-autonomous wave equations 
$$
u_{t t}+h\left(u_t\right)-\Delta u+f(u)=g(x, t),
$$
focusing on the existence of uniform attractors under critical growth conditions. They developed a novel methodology to demonstrate uniform asymptotic compactness, a key property for attractor existence. Yang \cite{LY13} subsequently extended this approach to non-autonomous plate equations with critical nonlinearity
$$
u_{t t}+a(x) u_t+\Delta^2 u+\lambda u+f(u)=g(x, t).
$$
More recently, Wang, Li, and Qin \cite{YW28} turned their attention to the upper-semicontinuity of uniform attractors for non-autonomous nonclassical diffusion equations with critical nonlinearity
$$
\partial_t u-\varepsilon \Delta \partial_t u-\Delta u+f(u)=g(x, t),
$$
their work is particularly noteworthy as it addresses the upper-semicontinuity of attractors, a property that ensures the stability of the attractor under small perturbations of the system parameters.

Our current study investigates a non-autonomous viscoelastic Kirchhoff plate equation with memory effects 	
$$
u_{t t}-\Delta u_{t t}+a_\epsilon(t) u_t+\alpha \Delta^2 u-\int_0^{\infty} \mu(s) \Delta^2 u(t-s) \mathrm{d} s-\Delta u_t+f(u)=g(x,t).
$$
From a physical perspective, this model describes the transverse vibrations of a thick plate composed of viscoelastic material. The term $-\Delta u_{t t}$ represents the rotational inertia, which is essential when the thickness of the plate is not negligible. The integral memory term captures the hereditary characteristics of the material-common in polymers or biological tissues-where the stress depends on the entire history of the strain, rather than just the instantaneous state. The interplay between the strong damping $-\Delta u_t$ and the memory term further characterizes the complex dissipation mechanism of the system.

Mathematically, the presence of the time-dependent external force $g(x, t)$ and damping coefficient $a_\epsilon(t)$ renders the classical theory of global attractors inapplicable. This necessitates the study of the uniform attractor, which captures the minimal closed set attracting all trajectories uniformly with respect to the initial time. Establishing the uniformly absorbing set and verifying uniform asymptotic compactness are crucial steps to guarantee the existence of such an attractor. Furthermore, we explore the upper semicontinuity of the family of uniform attractors as the perturbation parameter $\epsilon \rightarrow 0^{+}$. This investigation essentially addresses the robustness (or structural stability) of the system, ensuring that the long-term dynamic structures do not collapse under small fluctuations of the damping parameters, which is of significant importance in practical engineering control and optimization.

This work builds on the methodologies established by Sun et al. \cite{CS6} and Yang \cite{LY13}, seeking to contribute to the understanding of the long-term behavior of non-autonomous Kirchhoff-type equations with complex interactions between damping, memory effects, and external forcing.

\begin{Remark}
	Before proceeding, it is worth clarifying the specific roles and necessity of the strong damping term $-\Delta u_t$ and the memory kernel $\mu$ in our model:
	
	(i) Role of $-\Delta u_t$ : The term $-\Delta u_t$ provides "strong damping" to the plate system, offering a crucial regularizing effect. While the weak damping $a_\epsilon(t) u_t$ and the viscoelastic memory term also contribute to dissipation, the presence of $-\Delta u_t$ significantly facilitates the proof of uniform asymptotic compactness. If this term were absent, the regularity of the uniform attractor would likely be lower, and establishing its existence would require more restrictive growth conditions on the nonlinear term $f(u)$ (e.g., subcritical growth only).
	
	(ii) Role of memory $\mu$ : The memory term is central to the novelty of this work. If $\mu \equiv 0$, system (1.1) reduces to a standard non-autonomous Kirchhoff plate equation without hereditary effects. The inclusion of the memory term necessitates the introduction of a new history variable $\eta^t$, which fundamentally changes the phase space from $V_2 \times V_1$ to the infinitedimensional space $\mathcal{H}=V_2 \times V_1 \times \mathcal{M}_2$. This introduces complex infinite-dimensional dynamics and analytical challenges that are not present in memory-free systems.
\end{Remark}

Adhering to the framework outlined in \cite{CD3,CG1,VP2}, we introduce a new variable $\eta^t$ to the system, representing the history of relative displacement. The definition is as follows
\begin{equation}\label{eq1.4}
	\eta^t(x, s)=u(x, t)-u(x, t-s), \quad(x, s) \in \Omega \times \mathbb{R}^{+}, t \geq 0.
\end{equation}

Differentiating (\ref{eq1.4}) with respect to $t$, we have
$$
\eta_t^t(x, s)=-\eta_s^t(x, s)+u_t(x, t), \quad(x, s) \in \Omega \times \mathbb{R}^{+}, t \geq 0,
$$
and we can take initial condition $(t=0)$ as
$$
\eta^0(x, s)=u_0(x, 0)-u_0(x,-s), \quad(x, s) \in \Omega \times \mathbb{R}^{+}.
$$

Thus, the original memory term can be rewritten as
$$
\int_0^{\infty} \mu(s) \Delta^2 u(t-s) \mathrm{d} s=\left(\int_0^{\infty} \mu(s) \mathrm{d} s\right) \Delta^2 u(t)-\int_0^{\infty} \mu(s) \Delta^2 \eta^t(s) \mathrm{d} s.
$$

To streamline the analysis, we set $\alpha := 1 + \int_0^{\infty} \mu(s),\mathrm{d}s$, under which equation (\ref{eq1.1}) can be reformulated as
\begin{equation}\label{eq1.5}
	\left\{\begin{array}{l}
		u_{t t}-\Delta u_{t t}+a_\epsilon(t) u_t+\Delta^2 u+\int_0^{\infty} \mu(s) \Delta^2 \eta^t(s) \mathrm{d} s-\Delta u_t+f(u)=g(x,t), \\
		\eta_t^t(x, s)=-\eta_s^t(x, s)+u_t(x, s),
	\end{array}\right.
\end{equation}
with boundary condition
\begin{equation}\label{eq1.6}
	u=\Delta u=0 \quad \text { on } \quad \partial \Omega \times \mathbb{R}^{+}, \quad \eta^t=\Delta \eta^t=0 \quad \text { on } \quad \partial \Omega \times \mathbb{R}^{+},
\end{equation}
and initial conditions
\begin{equation}\label{eq1.7}
	u(x, 0)=u_0(x), \quad u_t(x, 0)=u_1(x), \quad \eta^t(x, 0)=0, \quad \eta^0(x, s)=\eta_0(x, s),
\end{equation}
where
$$
	\begin{cases}
		u_0(x)=u_0(x, 0), & x \in \Omega, \\ u_1(x)=\left.\partial_t u_0(x, t)\right|_{t=0}, & x \in \Omega, \\ \eta_0(x, s)=u_0(x, 0)-u_0(x,-s), & (x, s) \in \Omega \times \mathbb{R}^{+}.
	\end{cases}
$$

This paper aims to establish the existence of solutions for the non-autonomous problem described by equations (\ref{eq1.1})--(\ref{eq1.3}) and to investigate the upper semicontinuity of uniform attractors with respect to perturbations in the damping term $a_{\epsilon}(t)u_t$. Our analysis covers both subcritical and critical exponents for the function f, offering a comprehensive approach to the problem.

The main contributions of this paper are as follows:
 
 (i) We address the asymptotic dynamics of the non-autonomous viscoelastic Kirchhoff plate equation with memory effects. The challenge lies in handling the nonlinear term $f(u)$ and the memory term $-\int_0^{\infty} \mu(s) \Delta^2 u(t-s) \mathrm{d} s$. We employ a multiplier approach to construct a uniformly absorbing set. However, the complexity of these terms necessitates a meticulous analytical strategy.
 
 (ii) We verify the conditions for uniform asymptotic compactness using contractive functionals. Compared to existing literature, our study introduces several novel aspects. For instance, while previous works have primarily focused on autonomous Kirchhoff equations or non-autonomous equations without memory effects, our research incorporates both non-autonomous forcing and memory terms. This combination presents unique challenges, particularly in handling the nonlinear term $f(u)$ and the memory term $-\int_0^{\infty} \mu(s) \Delta^2 u(t-s) \mathrm{d} s$. We address these challenges by employing a multiplier approach to construct a uniformly absorbing set and by verifying the conditions for uniform asymptotic compactness using contractive functionals. Our methodology, which includes the application of a quasi-stability inequality for subcritical growth and a novel approach for critical growth, extends the techniques established by Chueshov and Lasiecka \cite{IC29, IC30} and Sun, Cao, and Duan \cite{CS6}.
 
 (iii) We address the open question of upper semicontinuity of uniform attractors in the presence of arbitrary time-dependent perturbations. We derive the upper semicontinuity of uniform attractors as the damping terms tend to zero, considering cases where $a_{\epsilon}$ is dependent solely on $\epsilon$ or where $a_\epsilon(t)$ is bounded. To the best of our knowledge, this is the first study to tackle the non-autonomous plate equation with critical nonlinearity, marking a significant advancement in the theoretical understanding of such systems.
 
 The paper is structured as follows. In Section 2, we present the functional setting, assumptions, and preliminary tools. Section 3 establishes the existence of a uniformly (w.r.t. $\sigma \in \Sigma$) absorbing set, as in Theorem \ref{th3.1}. Section 4 is devoted to proving the uniformly (w.r.t. $\sigma \in \Sigma$) asymptotic compactness of the solution process, see Theorem \ref{th4.4}. Section 5 deals with the upper semicontinuity of the family of uniform attractors with respect to small perturbations in the damping term, as in Theorem \ref{th5.1}.

\section{Preliminaries}

 \quad In this section, we initiate by presenting several function spaces, drawing parallels with those detailed in \cite{XP7}. Following this, essential assumptions are outlined to frame our analysis. Lastly, we revisit fundamental concepts and key results from the uniform attractor theory, as documented in \cite{VV4,IM5,CS6}, setting the stage for our subsequent discussions.

\textbf{\large {2.1 Some spaces}}

We start this section with introducing the following Hilbert spaces
$$
V_0=L^2(\Omega), \quad V_1=H_0^1(\Omega), \quad V_2=H^2(\Omega) \cap H_0^1(\Omega),
$$
and
$$
V_3=\left\{u \in H^3(\Omega) \mid u=\Delta u=0 \text { on } \partial \Omega\right\},
$$
equipped with respective inner products and norms,
$$
\begin{aligned}
	& (u, v)_{V_1}=(\nabla u, \nabla v) \quad \text { and } \quad\|u\|_{V_1}=\|\nabla u\| \\
	& (u, v)_{V_2}=(\Delta u, \Delta v) \quad \text { and } \quad\|u\|_{V_2}=\|\Delta u\| \\
	& (u, v)_{V_3}=(\nabla \Delta u, \nabla \Delta v) \quad \text { and } \quad\|u\|_{V_3}=\|\nabla \Delta u\| .
\end{aligned}
$$

As usual, $\|\cdot\|_2$ denotes the $L^2$-norms, for simplify, we denote $\|\cdot\|_2$ as $\|\cdot\|$, the notation $(\cdot, \cdot)$ is used to represent either the standard inner product in $L^2$ or the duality pairing between a Banach space $V$ and its topological dual $V'$. The constants $\lambda_0, \lambda_1, \lambda_2>0$ represent the embedding constants
\begin{equation}\label{eq2.1}
	\lambda_0\|u\|_2^2 \leq\|\nabla u\|_2^2, \quad \lambda_1\|u\|_2^2 \leq\|\Delta u\|_2^2, \quad \lambda_2\|\nabla u\|_2^2 \leq\|\Delta u\|_2^2 \quad \text { for } u \in V_2 .
\end{equation}

To regard the relative displacement $\eta^t$ as a new variable, the weighted $L^2$-spaces are defined
$$
\mathcal{M}_i:=L_\mu^2\left(\mathbb{R}^{+} ; V_i\right)=\left\{\xi: \mathbb{R}^{+} \rightarrow V_i \mid \int_0^{\infty} \mu(s)\|\xi(s)\|_{V_i}^2 \mathrm{~d} s<\infty\right\}, \quad i=0,1,2,3,
$$
which are Hilbert spaces endowed with inner products and norms
$$
\left(\xi_1, \xi_2\right)_{\mu, i}=\int_0^{\infty} \mu(s)\left(\xi_1(s), \xi_2(s)\right)_{V_i} \mathrm{~d} s,
$$
and
$$
\|\xi\|_{\mu, i}^2=\int_0^{\infty} \mu(s)\|\xi(s)\|_{V_i}^2 \mathrm{~d} s, \quad i=0,1,2,3,
$$
respectively.

Now let us introduce the phase spaces
\begin{equation}\label{eq2.2}
	\mathcal{H}=V_2 \times V_1 \times \mathcal{M}_2 \quad \text { and } \quad \mathcal{H}_1=V_3 \times V_2 \times \mathcal{M}_3,
\end{equation}
equipped with norms
$$
\|(u, v, \xi)\|^2_{\mathcal{H}}=\|\Delta u\|_2^2+\|\nabla v\|_2^2+\|\xi\|_{\mu, 2}^2,
$$
and
$$
\|(u, v, \xi)\|^2_{\mathcal{H}_1}=\|\nabla \Delta u\|_2^2+\|\Delta v\|_2^2+\|\xi\|_{\mu, 3}^2,
$$
respectively.

\textbf{\large {2.2 Assumptions}}

We now specify a set of assumptions concerning the nonlinear term $f$, the time-dependent coefficient $a_{\epsilon}(t)$, the source function $g$, and the memory kernel $\mu$.

Assume that $f \in C^1(\mathbb{R})$ satisfies
\begin{equation}\label{eq2.3}
	\left|f^{\prime}(s)\right| \leq M_f\left(1+|s|^p\right), \quad \forall s \in \mathbb{R},
\end{equation}
in addition, we assume 
\begin{equation}\label{eq2.4}
	-C_f \leq F(s) \leq f(s) s, \quad \forall s \in \mathbb{R},
\end{equation}
where $M_f, C_f$ are positive constants, $F(z)=\int_0^z f(s) d s$ with
\begin{equation}\label{eq2.5}
	0<p \leq \frac{4}{N-4}=p^* \quad \text{if} \quad N \geq 5 \quad\text{and} \quad 0<p<\infty \quad \text{if}\quad 1 \leq N \leq 4.
\end{equation}
We note that the condition (\ref{eq2.5}) provides $H^2(\Omega) \hookrightarrow L^{2(p+1)}(\Omega)$.

Inspired from \cite{MA27}, we assume $a_\epsilon: \mathbb{R} \rightarrow \mathbb{R}$ is a continuous and bounded function that satisfies
\begin{equation}\label{eq2.6}
	0<a_0 \leq a_\epsilon(t) \leq a_1, \forall \epsilon \in[0,1], \forall t \in \mathbb{R},
\end{equation}
for some positive constants $a_0$ and $a_1$.
For instance, the following is an example of functions depending on $\epsilon$ that satisfies the above conditions
\begin{equation}\label{eq2.7}
	a_\epsilon(t)=\frac{(1-\epsilon) t^2}{1+t^2}+\zeta, \quad \forall \epsilon \in[0,1], \quad \forall t \in \mathbb{R}, \quad \forall \zeta>0.
\end{equation}

Regarding the time-dependent external force term $g(x,t)$, we make the following assumption
\begin{equation}\label{eq2.8}
	g \in L_{w, l o c}^2\left(\mathbb{R} ; L^2(\Omega)\right),
\end{equation}
we denote by $L^2_{w,\mathrm{loc}}(\mathbb{R}; L^2(\Omega))$ the space $L^2_{\mathrm{loc}}(\mathbb{R}; L^2(\Omega))$ equipped with the topology of local weak convergence. For any function $g_0 \in L^2_b(\mathbb{R}; L^2(\Omega))$—i.e., $g_0$ is bounded under translations in the weak local topology—we define the associated set of temporal shifts by
\begin{equation}\label{eq2.9}
	\Sigma_0:=\left\{g_0(\cdot, \cdot+h) \mid h \in \mathbb{R}\right\} .
\end{equation}
and denote by $\Sigma$ the closure of $\Sigma_0$ in $L^2_{w,\mathrm{loc}}(\mathbb{R}; L^2(\Omega))$. It follows that every $g \in \Sigma$ remains in $L^2_b(\mathbb{R}; L^2(\Omega))$, meaning
\begin{equation}\label{eq2.10}
	\|g\|_{L^2_b}^2 := \sup_{t \in \mathbb{R}} \int_t^{t+1} \|g(s)\|_{L^2(\Omega)}^2 \, ds < \infty,
\end{equation}
where $\|\cdot\|_{L^2_b}$ denotes the norm in the space of translation-bounded functions $L^2_b(\mathbb{R}; L^2(\Omega))$.
\begin{Remark}
	To verify the translation compactness assumption and facilitate the understanding of the symbol space $\Sigma$, we provide a concrete example of the external force $g_0(x, t)$ :
	
	 Periodic case: Let $g_0(x, t)=h(x) \sin t$ where $h \in L^2(\Omega)$. Since $\sin t$ is periodic, the set of all its time-translations $\left\{g_0(\cdot, t+h) \mid h \in \mathbb{R}\right\}$ forms a closed curve (a circle) in the phase space, which is clearly compact.
\end{Remark} 

We next formulate the conditions imposed on the memory kernel $\mu$.
\begin{equation}\label{eq2.11}
	\mu \in C^1\left(\mathbb{R}^{+}\right) \cap L^1\left(\mathbb{R}^{+}\right), \quad \mu^{\prime}(s) \leq 0, \quad \mu(s) \geq 0,
\end{equation}
and there exist constants $\mu_0, \delta>0$ such that
\begin{equation}\label{eq2.12}
	\int_0^{\infty} \mu(s) \mathrm{d} s=\mu_0, \quad 0<\mu(0)<\infty ,
\end{equation}
and
\begin{equation}\label{eq2.13}
	\mu^{\prime}(s)+\delta \mu(s) \leq 0, \quad \forall s \in \mathbb{R}^{+}.
\end{equation}

\textbf{\large {2.3 Some basic results}}
\begin{Definition}\rm(\cite{VV4})\label{de2.1}
	 (i) A family of sets $\left\{U_\sigma(t, \tau) \mid t \geq \tau, \tau \in \mathbb{R}\right\}, \sigma \in \Sigma$ (parameter set) is said to be a family of processes acting on metric space $M$ if for each $\sigma \in \Sigma$, $\left\{U_\sigma(t, \tau) \mid t \geq \tau, \tau \in \mathbb{R}\right\}$ is a process acting on $M$, i.e., the two-parameter mappings from $M$ to $M$ satisfying
	 \begin{equation}\label{eq2.14}
	 	U_\sigma(t, s) U_\sigma(s, \tau)=U_\sigma(t, \tau), \quad \forall t \geq s \geq \tau, \tau \in \mathbb{R},
	 \end{equation}
	 \begin{equation}\label{eq2.15}
	  	U_\sigma(\tau, \tau)=I \text { (identity operator), } \quad \tau \in \mathbb{R}.
	  \end{equation}
	And the set $\Sigma$ is said to be the symbol space and $\sigma \in \Sigma$ to be a symbol.
	
	(ii) Let $\{T(t)\}_{t \geq 0}$ be a translation semigroup acting on $\Sigma$. The family of processes $\left\{U_\sigma(t, \tau)\right\}, \sigma \in \Sigma$ is said to be satisfying the translation identity if
	\begin{equation}\label{eq2.16}
			U_\sigma(t+s, \tau+s)=U_{T(s) \sigma}(t, \tau), \quad \forall \sigma \in \Sigma, t \geq \tau, \tau \in \mathbb{R}, s \geq 0.
	\end{equation}

	(iii) A bounded subset $B_0 \subset M$ is said to be a bounded uniformly (w.r.t. $\sigma \in \Sigma$ ) absorbing set of the family of processes $\left\{U_\sigma(t, \tau)\right\}, \sigma \in \Sigma$ if for any $\tau \in \mathbb{R}$ and bounded subset $B \subset M$ there exists a $T_0=T_0(B, \tau) \geq \tau$ such that
	\begin{equation}\label{eq2.17}
			\bigcup_{\sigma \in \Sigma} U_\sigma(t, \tau) B \subset B_0, \quad \forall t \geq T_0.
	\end{equation}
\end{Definition}
\begin{Definition}\rm(\cite{VV4})\label{de2.2}
	(i) We say that the family of processes $\left\{U_\sigma(t, \tau): \sigma \in \Sigma\right\}$, is uniformly (with respect to (w.r.t.) $\sigma \in \Sigma$ ) bounded if for any $B \in \mathcal{B}(E)$, the set
	$$
	\bigcup_{\sigma \in \Sigma} \bigcup_{\tau \in \mathbb{R}} \bigcup_{t \geq \tau} U_\sigma(t, \tau) B \in \mathcal{B}(E)
	$$
	(ii) A bounded set $B_0 \in \mathcal{B}(E)$ is said to be a bounded uniformly (w.r.t $\sigma \in \Sigma$ ) absorbing set for $\left\{U_\sigma(t, \tau): \sigma \in \Sigma\right\}\left(t \geq \tau, \tau \in \mathbb{R}^{+}\right)$if for any $\tau \in \mathbb{R}^{+}$and $B \in \mathcal{B}(E)$, there exists a time $T_0=T_0(B, \tau) \geq \tau$ such that
	$$
	\bigcup_{\sigma \in \Sigma} U_\sigma(t, \tau) B \subseteq B_0, \quad \text { for all } t \geq T_0.
	$$
\end{Definition}

Remember the definition of Hausdorff semi-distance between two subsets $A$ and $B$ of a metric space $(E, d)$ :
$$
\operatorname{dist}_E(A, B)=\sup _{a \in A} \inf _{b \in B} d(a, b)
$$
\begin{Definition}\rm(\cite{VV4})\label{de2.3}
	A closed set $\mathcal{A}_{\Sigma} \subset M$ is said to be the uniform (w.r.t. $\sigma \in \Sigma$ ) attractor of the family of processes $\left\{U_\sigma(t, \tau)\right\}, \sigma \in \Sigma$ if
	
	(i) (Uniform attractiveness) $\mathcal{A}_{\Sigma}$ uniformly (w.r.t. $\sigma \in \Sigma$ ) attracts all the bounded subsets in $M$, i.e., for every bounded subset $B \subset M$ and $\tau \in \mathbb{R}$,
	$$
	\lim _{t \rightarrow \infty} \sup _{\sigma \in \Sigma} \operatorname{dist}_M\left\{U_\sigma(t, \tau) B, \mathcal{A}_{\Sigma}\right\}=0,
	$$
	here, $\operatorname{dist}_M\{\cdot, \cdot\}$ is the Hausdorff semi-distance in $M$, i.e.,
	$$
	\operatorname{dist}_M\{A, B\}=\sup _{x \in A} \inf _{y \in B} d_M(x, y), \quad A, B \subset M.
	$$
	
	(ii) (Minimality) For any closed set $\mathcal{A}^{\prime} \subset M$, if $\mathcal{A}^{\prime}$ is of property (i), then $\mathcal{A}_{\Sigma} \subset \mathcal{A}^{\prime}$.
\end{Definition} 

\begin{Definition}\rm(\cite{VV4})\label{de2.4}
	 A process $\left\{U_\sigma(t, \tau)\right\}$ in a metric space $E$ is said to be is uniformly (w.r.t. $\sigma \in \Sigma$ ) asymptotically compact if, for each fixed $\tau \in \mathbb{R}$, each sequence $\left\{t_n\right\} \subset[\tau, \infty)$ with $t_n \xrightarrow{n \rightarrow \infty} \infty$, and each bounded sequence $\left\{u_n\right\} \subset E,\left\{\sigma_n\right\} \subset \Sigma$ the sequence $\left\{U_{\sigma_n}\left(t_n, \tau\right) u_n\right\}$ has a convergent subsequence in $E$.
\end{Definition}

Now, we begin by defining uniformly contractive functions and discussing their relationship with uniform asymptotic compactness, refering to \cite{IM5,CS6} for detailed insights. We then introduce the necessary function spaces, as in \cite{XP7}, and outline the assumptions that guide our analysis. Finally, we revisit key concepts and results from uniform attractor theory, drawing from \cite{VV4,IM5,CS6}, to establish a solid foundation for our subsequent discussions.

\begin{Definition}\rm(\cite{IM5,CS6})\label{de2.5}
	Let $E$ be a Banach space and $B$ be a bounded subset of $E, \Sigma$ be a symbol (or parameter) space. We call a function $\phi(\cdot, \cdot ; \cdot, \cdot)$, defined on $(E \times E) \times(\Sigma \times$ $\Sigma)$ to be a contractive function on $B \times B$ if for any sequence $\left\{x_n\right\}_{n=1}^{\infty} \subseteq B$ and any $\left\{\sigma_n\right\} \subseteq \Sigma$, there are subsequences $\left\{x_{n_k}\right\}_{k=1}^{\infty} \subset\left\{x_n\right\}_{n=1}^{\infty}$ and $\left\{\sigma_{n_k}\right\}_{k=1}^{\infty} \subset\left\{\sigma_n\right\}_{n=1}^{\infty}$ such that
	$$
	\lim _{k \rightarrow \infty} \lim _{l \rightarrow \infty} \phi\left(x_{n_k}, x_{n_l} ; \sigma_{n_k}, \sigma_{n_l}\right)=0.
	$$
\end{Definition}

We denote the collection of all contractive functions on 
$B \times B$ by $\operatorname{Contr}(B, \Sigma)$.
The subsequent lemma represents an adaptation of Khanmamedov's result \cite{AK9} for semi-processes, which is pivotal for establishing the uniform asymptotic compactness of the semi-process. For the proof, readers are referred to the works of Chueshov and Lasiecka \cite{IC10}, as well as Sun et al. \cite{CS6,CS11}.
 
 \begin{Lemma}\rm(\cite{CS6,CS11})\label{le2.6}
 	Let $\left\{U_\sigma(t, \tau)\right\}\left(t \geq \tau, \tau \in \mathbb{R}^{+}, \sigma \in \Sigma\right)$ be a process satisfying the translation identities (\ref{eq2.14}) and (\ref{eq2.15}) on Banach space $E$ and suppose that $\left\{U_\sigma(t, \tau)\right\}$ has a bounded uniformly (w.r.t. $\sigma \in \Sigma$ ) absorbing set $B_0 \subseteq E$. Moreover, suppose that for any $\varepsilon>0$, there exist $T=T\left(B_0, \varepsilon\right)>0$ and $\phi_T \in \operatorname{Contr}\left(B_0, \Sigma\right)$ such that
 	$$
 	\left\|U_{\sigma_1}(T, 0) x-U_{\sigma_2}(T, 0) y\right\| \leq \varepsilon+\phi_T\left(x, y ; \sigma_1, \sigma_2\right), \forall \sigma \in \Sigma, t \geq \tau, \tau \in \mathbb{R}^{+}.
 	$$
 	Then $\left\{U_\sigma(t, \tau)\right\}\left(t \geq \tau, \tau \in \mathbb{R}^{+}, \sigma \in \Sigma\right)$ is uniformly (w.r.t. $\sigma \in \Sigma$ ) asymptotically compact in $E$.	
 \end{Lemma}
  
  \begin{Remark}\rm(\cite{CS6,CS11})
  	(i) If there exists a uniformly contractive function for the family of processes $\left\{U_\sigma(t, \tau)\right\}$, then the family of processes is uniformly asymptotically compact.
  	
  	(ii) If $E$ is a uniformly convex Banach space, then the uniformly asymptotic compactness is equivalent to the family of processes possesses a contractive function.
  \end{Remark}
 
 \begin{Definition}\rm(\cite{CS6,CS11})\label{de2.7}
 	The set
 	$$
 	\omega_{\tau, \Sigma}(B)=\bigcap_{t \geq \tau} \overline{\bigcup_{\sigma \in \Sigma} \bigcup_{s \geq t} U_\sigma(s, \tau) B},
 	$$
 	where $B$ is a bounded subset of $E$ is said to be the uniform (w.r.t. $\sigma \in \Sigma$ ) $\omega$-limit set of $B$.
 \end{Definition}
 \begin{Lemma}\rm(\cite{CS6,CS11})\label{le2.8}
 	Let $\left\{U_\sigma(t, \tau): \sigma \in \Sigma\right\}$, be a family of uniform (w.r.t. $\sigma \in \Sigma$ ) asymptotically compact processes on a complete metric space $E$, then for any $\tau \in \mathbb{R}$ and any nonempty set $B \in \mathcal{E}$ we have
 	$$
 	\omega_{\tau, \Sigma}(B)=\omega_{0, \Sigma}(B),
 	$$
 	that is, $\omega_{\tau, \Sigma}(B)$ is independent of $\tau \in \mathbb{R}$.
 \end{Lemma}
 
 The next result relates uniform attractor with omega limit set. For its proof, see \cite{CS6}.
 \begin{Theorem}\rm(\cite{CS6,CS11})\label{th2.9}
 	Let $E$ be a complete metric space and $\left\{U_\sigma(t, \tau): \sigma \in \Sigma\right\}$ be a family of processes on $E$ which satisfies the translation identity (\ref{eq2.14})-(\ref{eq2.15}). Then, $\left\{U_\sigma(t, \tau): \sigma \in \Sigma\right\}$, has a compact uniform (w.r.t. $\sigma \in \Sigma$ ) attractor $A_{\Sigma}$ in $E$ and satisfies
 	$$
 	A_{\Sigma}=\omega_{0, \Sigma}\left(B_0\right)=\omega_{\tau, \Sigma}\left(B_0\right)=\bigcup_{B \in \mathcal{B}(E)} \omega_{\tau, \Sigma}(B), \quad \text { for all } \tau \in \mathbb{R},
 	$$
 	if and only if $\left\{U_\sigma(t, \tau): \sigma \in \Sigma\right\}$
 	
 	(i) has a bounded uniformly (w.r.t. $\sigma \in \Sigma$ ) absorbing set $B_0$; and
 	
 	(ii) is uniformly (w.r.t. $\sigma \in \Sigma$ ) asymptotically compact.
 \end{Theorem}
 \begin{Definition}\rm(\cite{AC12})\label{de2.10}
 	Let $E$ be a Banach space and $\epsilon$ be a parameter varying in the interval $[0,1]$. For each $\epsilon \in[0,1]$, let $\mathcal{A}_\epsilon$ be a subset of $E$. We say that the family $\left\{\mathcal{A}_\epsilon\right\}_{\epsilon \in[0,1]}$ is upper-semicontinuity at $\epsilon=0$ if
 	$$
 	\lim _{\epsilon \rightarrow 0^{+}} \operatorname{dist}_E\left(\mathcal{A}_\epsilon, \mathcal{A}_0\right)=0.
 	$$
 \end{Definition}
  
Similar to the autonomous cases (refer to \cite{MA27,XP7}), the existence and uniqueness results can be readily derived. The proof relies on the Galerkin approximation method (see \cite{CS6}), with time-dependent terms introducing no fundamental complexities.
\begin{Theorem}\rm(\cite{XP7})\label{th2.11}
	Assume that assumptions (\ref{eq2.3})-(\ref{eq2.13}) hold and consider $h \in V_0$. Then we have
	
	(i) If initial data $\left(u_0, u_1, \eta_0\right) \in \mathcal{H}$, then problem (\ref{eq1.5})-(\ref{eq1.7}) has a weak solution
	$$
	\left(u, u_t, \eta^t\right) \in C([0, T], \mathcal{H}), \quad \forall T>0,
	$$
	satisfying
	$$
	\begin{aligned}
		& u \in L^{\infty}\left(0, T ; V_2\right), \quad u_t \in L^{\infty}\left(0, T ; V_1\right), \\
		& (\mathrm{I}-\Delta) u_{t t} \in L_{l o c}^{\infty}\left(\mathbb{R}^{+}, V_2^{\prime}\right), \quad \eta^t \in L^{\infty}\left(0, T ; \mathcal{M}_2\right).
	\end{aligned}
	$$
	
	(ii) If initial data $\left(u_0, u_1, \eta_0\right) \in \mathcal{H}_1$, then problem (\ref{eq1.5})-(\ref{eq1.7}) has a stronger solution satisfying
	$$
	\begin{aligned}
		& u \in L^{\infty}\left(0, T ; V_3\right), \quad u_t \in L^{\infty}\left(0, T ; V_2\right), \\
		& (\mathrm{I}-\Delta) u_{t t} \in L_{l o c}^{\infty}\left(\mathbb{R}^{+}, V_1^{\prime}\right), \quad \eta^t \in L^{\infty}\left(0, T ; \mathcal{M}_3\right)
	\end{aligned}
	$$
	
	(iii) Let $z_1(t)=\left(u, u_t, \eta^t\right), z_2(t)=\left(v, v_t, \xi^t\right)$ be weak solutions of problem (\ref{eq1.5})-(\ref{eq1.7}) corresponding to initial data $z_1(0)=\left(u_0, u_1, \eta_0\right), z_2(0)=\left(v_0, v_1, \xi_0\right)$. Then one has
	$$
	\left\|z_1(t)-z_2(t)\right\|_{\mathcal{H}}^2 \leq C_T\left\|z_1(0)-z_2(0)\right\|_{\mathcal{H}}^2, \quad t \geq 0,
	$$
	for some constant $C_T=C\left(\left\|z_1(0)\right\|_{\mathcal{H}},\left\|z_2(0)\right\|_{\mathcal{H}}, T\right)>0$. In particular, problem (\ref{eq1.5})-(\ref{eq1.7}) has a unique weak solution.
\end{Theorem}

\section{Uniformly (w.r.t. $\sigma \in \Sigma$) absorbing set in $\mathcal{H}$}

 \quad In this section, we show the existence of uniformly absorbing set in $\mathcal{H}$. We first use the notations as in Chepyzhov and Vishik \cite{VV4}. Let $y(t)=\left(u(t), u_t(t), \eta^t(s)\right), y_\tau=\left(u_\tau^0, u_\tau^1,\eta^0\right), \mathcal{H}=\left(H^2(\Omega) \cap H_0^1(\Omega)\right) \times H_0^1(\Omega)\times \mathcal{M}_2$ with finite energy norm
$$
\|y\|^2_{\mathcal{H}}=\|\Delta u\|^2+\left\|\nabla u_t\right\|^2+\|\eta\|^2_{\mu,2} .
$$
Then system (\ref{eq1.5})--(\ref{eq1.7}) is equivalent to
\begin{equation}\label{eq3.1}
	\left\{\begin{array}{l}
		\partial_t u_t=\Delta u_{t t}-a_\epsilon(t) u_t-\Delta^2 u-\int_0^{\infty} \mu(s) \Delta^2 \eta^t(s) \mathrm{d} s+\Delta u_t-f(u)+g(x,t), \\
		\eta_t^t(x, s)=-\eta_s^t(x, s)+u_t(x, s),
	\end{array}\right.
\end{equation}
with boundary conditions
\begin{equation}\label{eq3.2}
	u=\Delta u=0 \quad \text { on } \quad \partial \Omega \times \mathbb{R}^{+}, \quad \eta^t=\Delta \eta^t=0 \quad \text { on } \quad \partial \Omega \times \mathbb{R}^{+},
\end{equation}
and initial conditions
\begin{equation}\label{eq3.3}
	u(x, 0)=u_0(x), \quad u_t(x, 0)=u_1(x), \quad \eta^t(x, 0)=0, \quad \eta^0(x, s)=\eta_0(x, s).
\end{equation}

Let

$A_{\sigma(t)}(y)=\left(\begin{array}{c}
	u_t \\
	(I-\Delta)^{-1}\left[-\left(a_\epsilon(t) I-\Delta\right) u_t-\alpha \Delta^2 u-\int_0^{\infty} \mu(s) \Delta^2 \eta^t(s) d s-f(u)+g(x, t)\right]\\
	-\eta_s+u_t
\end{array}\right).
$

 Then the nonautonomous system (\ref{eq3.1})
 can be rewritten in the operator form  
\begin{equation}\label{eq3.4}
	 \partial_t y=A_{\sigma(t)}(y),\left.\quad y\right|_{t=\tau}=y_\tau,
\end{equation}
where $\sigma(t)=g(x, t)$ is a symbol of Eq. (\ref{eq3.4}).

Next, we give the main result of this section as follows.
\begin{Theorem}\label{th3.1}
	Under assumptions (\ref{eq2.3})--(\ref{eq2.13}), the family of processes $\left\{U_\sigma(t, \tau)\right\}, \sigma \in \Sigma$, corresponding to (\ref{eq1.5})--(\ref{eq1.7}) has a bounded uniformly (w.r.t. $\sigma \in \Sigma$ ) absorbing set $B_0$ in $\mathcal{H}$.
	
	\rm\textbf{Proof.} Consider the energy of the system (\ref{eq1.5})--(\ref{eq1.7}) at time $t$, given by the functional $E: \mathcal{H} \rightarrow \mathbb{R}$
	\begin{equation}\label{eq3.5}
		E(t)=\frac{1}{2}\left\|u_t\right\|^2+\frac{1}{2}\left\|\nabla u_t\right\|^2+\frac{1}{2}\|\Delta u\|^2+\frac{1}{2}\left\|\eta^t\right\|_{\mu, 2}^2+\int_{\Omega}F(u(t))dx.
	\end{equation}
	
 Here, $g_0$ denotes the fixed external force (the generating function) given in the problem statement. From the definition of $\Sigma$, we know that for all $\sigma \in \Sigma$
 \begin{equation}\label{eq3.6}
 	\|\sigma\|_{L_b^2}^2 \leqslant\left\|g_0\right\|_{L_b^2}^2.
 \end{equation}

The time derivative of the energy functional $E(t)$ defined in (\ref{eq3.5}) is calculated by separating the terms
$$
	\frac{d}{dt}E(t)=\frac{d}{d t}\left(\frac{1}{2}\left\|u_t\right\|^2+\frac{1}{2}\left\|\nabla u_t\right\|^2\right)+\frac{d}{d t}\left(\frac{1}{2}\|\Delta u\|^2+\int_{\Omega} F(u) d x\right)+\frac{d}{d t}\left(\frac{1}{2}\left\|\eta^t\right\|_{\mu, 2}^2\right).
$$
Applying (\ref{eq1.5}) and integration by parts, we can obtain 
\begin{equation*}
	\begin{aligned}
			\frac{d}{d t}\left(\frac{1}{2}\left\|u_t\right\|^2+\frac{1}{2}\left\|\nabla u_t\right\|^2\right)&=\left(u_t, u_{t t}\right)+\left(\nabla u_t, \nabla u_{t t}\right)\\
			&=\left((I-\Delta) u_{t t}, u_t\right)\\
			&=\left(-a_\epsilon(t) u_t-\Delta^2 u-\int_0^{\infty} \mu(s) \Delta^2 \eta^t(s) \mathrm{d} s+\Delta u_t-f(u)+g(x, t), u_t\right),
	\end{aligned}
\end{equation*}
$$
\frac{d}{d t}\left(\frac{1}{2}\|\Delta u\|^2+\int_{\Omega} F(u) d x\right)=\left(\Delta^2 u, u_t\right)+\left(f(u), u_t\right),
$$
and
$$
\frac{d}{d t}\left(\frac{1}{2}\left\|\eta^t\right\|_{\mu, 2}^2\right)=-\left(\eta^t, \eta_{s}^t\right)_{\mu, 2}+\left(\int_0^{\infty} \mu(s) \Delta^2 \eta^t(s)ds, u_t\right).
$$

Summing the above three expressions, we directly obtain 
\begin{equation}\label{eq3.7}
	\frac{d}{dt}E(t)=-\left\|\nabla u_t\right\|^2-(\eta^t,\eta^t_s)_{\mu,2}-(a_{\epsilon}(t)u_t,u_t)+(g(t),u_t).
\end{equation}
Thanks to assumptions (\ref{eq2.11})--(\ref{eq2.13}), we obtain
\begin{equation}\label{eq3.8}
	\begin{aligned}
		\left(\eta^t, \eta_s^t\right)_{\mu, 2} & =\frac{1}{2} \int_0^{\infty} \mu(s) \frac{\mathrm{d}}{\mathrm{~d} t}\left\|\Delta \eta^t(s)\right\|_2^2 \mathrm{~d} s \\
		& =\left[\frac{1}{2} \mu(s)\left\|\Delta \eta^t(s)\right\|_2^2\right]_0^{\infty}-\frac{1}{2} \int_0^{\infty} \mu^{\prime}(s)\left\|\Delta \eta^t(s)\right\|_2^2 \mathrm{~d} s,
	\end{aligned}
\end{equation}
and
$$
  \lim _{s \rightarrow \infty} \mu(s)=0.
$$

According to the definition of $\eta^t(x, s)$, we can easily see that
$$
\eta^t(x, 0)=u(x, t)-u(x, t-0)=0,
$$
which implies that
$$
\left\|\Delta \eta^t(0)\right\|_2^2=0.
$$
Thus
$$
-\left(\eta^t, \eta_s^t\right)_{\mu, 2}=\frac{1}{2} \int_0^{\infty} \mu^{\prime}(s)\left\|\Delta \eta^t(s)\right\|_2^2 \mathrm{~d} s
$$
which, together with (\ref{eq2.13}), gives
\begin{equation}\label{eq3.9}
	-\left(\eta^t, \eta_s^t\right)_{\mu, 2} \leq -\frac{\delta}{2}\left\|\eta^t\right\|_{\mu, 2}^2
\end{equation}

From the assumption (\ref{eq2.1}) and (\ref{eq2.6}), we get
\begin{equation}\label{eq3.10}
	-(a_{\epsilon}(t)u_t,u_t)\leq -a_0\|u_t\|^2,
\end{equation}

Using the Young's inequality, we arrive at
\begin{equation}\label{eq3.11}
	(g(t),u_t)\leq \frac{1}{4a_0}\|g\|^2+a_0\|u_t\|^2.
\end{equation}

Combining (\ref{eq3.9})--(\ref{eq3.11}) with (\ref{eq3.7}), we obtain
\begin{equation}\label{eq3.12}
	\frac{d}{dt}E(t)\leq-\left\|\nabla u_t\right\|^2-\frac{\delta}{2}\left\|\eta^t\right\|_{\mu, 2}^2+\frac{1}{4a_0}\|g\|^2.
\end{equation}

Next, let us define
\begin{equation}\label{eq3.13}
	\Psi(t)=\int_{\Omega} u_t(t) u(t) \mathrm{d} x-\int_{\Omega} \triangle u_t(t) u(t) \mathrm{d} x.
\end{equation}
 
From (\ref{eq3.7}) and assumption (\ref{eq2.4}),
we have
\begin{equation}\label{eq3.14}
	E(t)+C_f|\Omega|\geq\frac{1}{2}\|\nabla u_t\|^2+\frac{1}{2}\|\Delta u\|^2+\frac{1}{2}\left\|\eta^t\right\|_{\mu, 2}^2=\frac{1}{2}\|(u,u_t,\eta^t)\|^2_{\mathcal{H}}.
\end{equation}

Considering the definition of $\Psi(t)$, we obtain
\begin{equation}\label{eq3.15}
	\begin{aligned}
		|\Psi(t)| & \leq \frac{1}{2}\left\|u_t\right\|^2+\frac{1}{2}\|u\|^2+\frac{1}{2}\left\|\nabla u_t\right\|^2+\frac{1}{2}\|\nabla u\|^2 \\
		& \leq \frac{1}{2}\left(\frac{1}{\lambda_0}+1\right)\left\|\nabla u_t\right\|^2+\frac{1}{2}\left(\frac{1}{\lambda_1}+\frac{1}{\lambda_2}\right)\|\Delta u\|^2 \\
		& \leq C_1\left(E(t)+C_f|\Omega|\right),
	\end{aligned}
\end{equation}
and
\begin{equation}\label{eq3.16}
	\begin{aligned}
		\frac{d}{dt}\Psi(t)&=\int_{\Omega}\left(u_{t t}-\Delta u_{t t}\right) u \mathrm{~d} x+\left\|u_t\right\|_2^2+\left\|\nabla u_t\right\|_2^2\\
		 & =-E(t)+\frac{3}{2}\|u_t\|^2+\frac{3}{2}\|\nabla u_t\|^2-\frac{1}{2}\|\Delta u\|^2+\frac{1}{2}\left\|\eta^t\right\|_{\mu, 2}^2+\int_{\Omega} F(u)-f(u)u \mathrm{~d} x\\
		 &\quad-\int_0^{\infty} \mu(s)\left(\Delta \eta^t(s), \Delta u(t)\right) \mathrm{d} s-\int_{\Omega}a_{\epsilon}(t)u_tudx -\int_{\Omega} \nabla u_t \cdot \nabla u \mathrm{~d} x+\int_{\Omega}g(t) u \mathrm{d} x.
	\end{aligned}
\end{equation}

From assumption (\ref{eq2.4}), we know
\begin{equation}\label{eq3.17}
	\int_{\Omega} F(u)-f(u)u \mathrm{~d} x\leq 0.
\end{equation}

 By the Young's inequality, given $\epsilon_1>0$, we get 
 \begin{equation}\label{eq3.18}
 	\begin{aligned}
 		\left|-\int_0^{\infty} \mu(s)\left(\Delta \eta^t(s), \Delta u(t)\right) \mathrm{d} s\right| & \leq \int_0^{\infty} \mu(s)\left(\frac{1}{4 \epsilon_1}\left\|\Delta \eta^t\right\|_2^2+\epsilon_1\|\Delta u\|_2^2\right) \mathrm{d} s \\
 		& =\epsilon_1\left(\int_0^{\infty} \mu(s) \mathrm{d} s\right)\|\Delta u\|_2^2+\frac{1}{4 \epsilon_1} \int_0^{\infty} \mu(s)\left\|\Delta \eta^t\right\|_2^2 \mathrm{~d} s \\
 		& =\epsilon_1\mu_0\|\Delta u\|_2^2+\frac{1}{4 \epsilon_1}\left\|\eta^t\right\|_{\mu, 2}^2.
 	\end{aligned}
 \end{equation}
	
By assumption (\ref{eq2.6}), we know
\begin{equation}\label{eq3.19}
	\left|\int_{\Omega}a_{\epsilon}(t)u_tudx\right|\leq a_1|(u_t,u)|\leq\frac{a_1^2\epsilon_1}{\lambda_1}\|\Delta u\|^2+\frac{1}{4\epsilon_1\lambda_0}\|\nabla u_t\|^2.
\end{equation}
 
 Using the Young's inequality, we get
 \begin{equation}\label{eq3.20}
 	\int_{\Omega}g(t)udx\leq\frac{1}{4\epsilon_1}\|g\|^2+\frac{\epsilon_1}{\lambda_1}\|\Delta u\|^2.
 \end{equation}

Inserting (\ref{eq3.17})--(\ref{eq3.20}) into (\ref{eq3.16}), we have
\begin{equation}\label{eq3.21}
	\begin{aligned}
		\frac{d}{dt}\Psi(t) \leq & -E(t)+\left(\frac{3}{2 \lambda_0}+\frac{3}{2}+\frac{1}{4\epsilon_1\lambda_0}+\frac{1}{4 \epsilon_1}\right)\left\|\nabla u_t\right\|^2+\left(\frac{1}{2}+\frac{1}{ 4\epsilon_1}\right)\left\|\eta^t\right\|_{\mu, 2}^2 \\
		& -\left[\frac{1}{2}-\epsilon_1\left(\mu_0+\frac{a_1^2}{\lambda_1}+\frac{1}{\lambda_1}+\frac{1}{\lambda_1}\right)\right]\|\Delta u\|^2.
	\end{aligned}
\end{equation}

 Choosing $\epsilon_1$ small enough such that
 $$
 	\frac{1}{2}-\epsilon_1\left(\mu_0+\frac{a_1^2}{\lambda_1}+\frac{1}{\lambda_1}+\frac{1}{\lambda_1}\right)>0,
 $$
 and denoting $C_2=\frac{3}{2 \lambda_0}+\frac{3}{2}+\frac{1}{4\epsilon_1\lambda_0}+\frac{1}{4 \epsilon_1}$, $C_3=\frac{1}{2}+\frac{1}{ 4\epsilon_1}$, we obtain
 \begin{equation}\label{eq3.22}
 	\frac{d}{dt}\Psi(t)\leq -E(t)+C_2\|\nabla u_t\|^2+C_3\left\|\eta^t\right\|_{\mu, 2}^2+\frac{1}{4\epsilon_1}\|g\|^2.
 \end{equation}

Then let us consider the perturbed energy
\begin{equation}\label{eq3.23}
	 E_{\varepsilon}(t)=E(t)+\varepsilon \Psi(t), \quad \varepsilon>0.
\end{equation}

We know from (\ref{eq3.15})
\begin{equation}\label{eq3.24}
	\left|E_{\varepsilon}(t)-E(t)\right| \leq \varepsilon C_1\left(E(t)+C_f|\Omega|\right), \quad \forall t \geq 0, \quad \forall \varepsilon>0.
\end{equation}	

We select $\varepsilon_1=\min \left\{\frac{1}{C_2}, \frac{\delta}{2 C_3}\right\}$, which is positive by our assumption that $\delta>0$. By combining equations (\ref{eq3.12}) and (\ref{eq3.22}), we can deduce that
\begin{equation}\label{eq3.25}
	\begin{aligned}
		\frac{d}{dt}E_{\varepsilon}(t) & =\frac{d}{dt}E(t)+\frac{d}{dt}\varepsilon \Psi(t) \\
		& \leq-\varepsilon E(t)-\left(1-\varepsilon C_2\right)\left\|\nabla u_t\right\|_2^2-\left(\frac{\delta}{2}-\varepsilon C_3\right)\left\|\eta^t\right\|_{\mu, 2}^2+\left(\frac{1}{4a_0}+\frac{1}{4\epsilon_1}\right)\|g\|^2 \\
		& \leq-\varepsilon E(t)+\left(\frac{1}{4a_0}+\frac{1}{4\epsilon_1}\right)\|g\|^2, \quad \varepsilon \in\left(0, \varepsilon_1\right] .
	\end{aligned}
\end{equation}

Denoting $C_4=\frac{1}{4a_0}+\frac{1}{4_{\epsilon_1}}$,
then
\begin{equation}\label{eq3.26}
	\frac{d}{dt}E_{\varepsilon}(t)\leq -\varepsilon E(t)+C_4\|g\|^2.
\end{equation}

  Define $\varepsilon_2 := \min\left\{\frac{1}{2C_1}, \varepsilon_1\right\}$. Then, for any $\varepsilon \in (0, \varepsilon_2]$, it follows from equation (\ref{eq3.24}) that
 \begin{equation}\label{eq3.27}
   \frac{1}{2} E(t)-\frac{1}{2}C_f|\Omega|\leq E_{\varepsilon}(t) \leq \frac{3}{2} E(t)+\frac{1}{2}C_f|\Omega|.
 \end{equation}

From (\ref{eq3.27}), we see that
\begin{equation}\label{eq3.28}
	E_{\varepsilon}^{\prime}(t) \leq-\frac{2 \varepsilon}{3} E_{\varepsilon}(t)+\frac{\varepsilon}{3}C_f|\Omega|+C_4\|g\|^2,
\end{equation}
 the following estimate holds for all $t \geq 0$ and $\varepsilon \in\left(0, \varepsilon_2\right]$ 
 \begin{equation}\label{eq3.29}
 	\begin{aligned}
 		E_{\varepsilon}(t) & \leq E_{\varepsilon}(0) \mathrm{e}^{-\frac{2 \varepsilon}{3} t}+\int_{0}^{t}\left(\frac{\varepsilon}{3}C_f|\Omega|+C_4\|g\|^2\right)\mathrm{e}^{-\frac{2 \varepsilon}{3} (t-s)}ds \\
 		&\leq\left[E_{\varepsilon}(0)-\frac{1}{2}\left(C_f|\Omega|+\frac{3C_4}{\varepsilon}\|g_0\|^2_{L_b^2}\right)\right] \mathrm{e}^{-\frac{2 \varepsilon}{3} t}+\frac{1}{2}\left(C_f|\Omega|+\frac{3C_4}{\varepsilon}\|g_0\|^2_{L_b^2}\right) .
 	\end{aligned}
 \end{equation}

Using the left-hand side of (\ref{eq3.27}) and (\ref{eq3.29}), we get
\begin{equation}\label{eq3.30}
	E(t)\leq 2E_{\varepsilon}(0)\mathrm{e}^{-\frac{2 \varepsilon}{3} t}+2C_f|\Omega|+\frac{3C_4}{\varepsilon}\|g_0\|^2_{L_b^2},
\end{equation}
from the right-hand side of (\ref{eq3.27}) and (\ref{eq3.30}), we have 
\begin{equation}\label{eq3.31}
	E(t)\leq 3E(0)\mathrm{e}^{-\frac{2 \varepsilon}{3} t}+3C_f|\Omega|+\frac{3C_4}{\varepsilon}\|g_0\|^2_{L_b^2}.
\end{equation}

 Combining (\ref{eq3.14}) with (\ref{eq3.31}), we  conclude that
 \begin{equation}\label{eq3.32}
 	\left\|\left(u(t), u_t(t), \eta^t\right)\right\|_{\mathcal{H}}^2 \leq C E(0) \mathrm{e}^{-\frac{2 \varepsilon}{3} t}+C\left(C_f|\Omega|+\|g_0\|^2_{L_b^2}\right),
 \end{equation}
where $C=\max\{8,\frac{6C_4}{\varepsilon}\}$.

Consequently, one arrives at a uniform absorbing set $B_0 \subset \mathcal{H}$, valid for all $\sigma \in \Sigma$,
$$
B_0=\left\{y=\left(u, u_t,\eta^t\right) \mid\|y\|_{\mathcal{H}}^2 \leqslant \rho_0^2\right\},
$$
where $\rho_0=2C\left(C_f|\Omega|+\|g_0\|^2_{L_b^2}\right)$. In other words, given any bounded set $B \subset \mathcal{H}$, there exists a time $t_0 = t_0(B) \geq 0$ such that
$$
\bigcup_{\varepsilon\in[0,1]}\bigcup_{\sigma \in \Sigma} U_{\sigma}^{\varepsilon}(t, 0) B \subset B_0, \quad \forall t \geqslant t_0.
$$
Then one has 
\begin{equation}\label{eq3.33}
		\left\|\left(u(t), u_t(t), \eta^t\right)\right\|_{\mathcal{H}} \leq C_B, \quad \forall t \geq 0,	
\end{equation}
where $C_B$ is a constant depending on $B$.
 The proof is complete.
\end{Theorem}
 $\hfill\qedsymbol$
 
\section{Uniformly (w.r.t. $\sigma \in \Sigma$) asymptotic compactness in $\mathcal{H}$}

\quad In this section, we initially derive certain a priori estimates concerning energy inequalities, adopting the approach outlined in \cite{AK14,LY13}. Subsequently, using Lemma \ref{le2.6}, we establish the uniform (w.r.t. $\sigma \in \Sigma)$ asymptotic compactness in $\mathcal{H}$.

\textbf{\large {4.1 A priori estimates}}

This subsection is primarily devoted to proving equation (\ref{eq4.2}), which plays a key role in verifying the uniform asymptotic compactness of the process with respect to $\sigma \in \Sigma$.

For each $i = 1, 2$, consider initial data $\left(u_0^i, u_1^i, \eta^{0i}\right) \in B_0$, and denote by $\left(u_i(t), u_{i_t}(t), \eta^{ti}\right)$ the solution associated with the symbol $\sigma_i$ and the prescribed initial condition. That is, $\left(u_i(t), u_{i_t}(t), \eta^{ti}\right)$ satisfies the following system
\begin{equation}\label{eq4.1}
	\left\{\begin{array}{l}
		u_{t t}-\Delta u_{tt}+a_{\epsilon}(t) u_t+\Delta^2 u+\int_0^{\infty} \mu(s) \Delta^2 \eta^t(s) \mathrm{d} s-\Delta u_t+f(u)=\sigma_i(x, t), \\
		\eta_t^t=-\eta_s^t+u_t,\\
			u=\Delta u=0 \quad \text { on } \quad \partial \Omega \times \mathbb{R}^{+}, \quad \eta^t=\Delta \eta^t=0 \quad \text { on } \quad \partial \Omega \times \mathbb{R}^{+},\\
		 	u(x, 0)=u_0(x), \quad u_t(x, 0)=u_1(x), \quad \eta^t(x, 0)=0, \quad \eta^0(x, s)=\eta_0(x, s).
	\end{array}\right.
\end{equation}
\begin{Lemma}\label{le4.1}
	Assume the condition (\ref{eq2.3})--(\ref{eq2.13}) is satisfied. Then for any fixed $T>0$, there exist a constant $C_T$ and a function $\phi_T=\phi_T\left(y_1,y_2; \sigma_1, \sigma_2\right)$, $y_1=(u_0^1,u_1^1,\eta^{01})$, $y_1=(u_0^2,u_1^2,\eta^{02})$, such that
	\begin{equation}\label{eq4.2}
		\left\|U_{\sigma_1}(T,0)y_1-U_{\sigma_2}(T,0)y_2\right\|_{\mathcal{H}} \leqslant \varepsilon_T+\phi_T\left(y_1,y_2; \sigma_1, \sigma_2\right),
	\end{equation}
	where $\varepsilon_T$ and $\phi_T$ depend on $T$.
	
	\rm\textbf{Proof.} For convenience, we denote
	$$
	g_i(t)=\sigma_i(x, t), \quad t \geqslant 0, i=1,2,
	$$
	and
	$$
	w(t)=u_1(t)-u_2(t), \quad\xi^t=\eta^{t1}-\eta^{t2}.
	$$
	
	Then $w(t)$ satisfies
	\begin{equation}\label{eq4.3}
			\left\{\begin{array}{l}
			w_{t t}-\Delta u_{tt}+a_{\epsilon}(t) w_t+\Delta^2 w+\int_0^{\infty} \mu(s) \Delta^2 \xi^t(s) \mathrm{d} s-\Delta u_t\\\quad+f\left(u_1(t)\right)-f\left(u_2(t)\right)=g_1(t)-g_2(t), \\
			\xi_t^t=-\xi_s^t+w_t,\\
			w=\Delta w=0 \quad \text { on } \quad \partial \Omega \times \mathbb{R}^{+}, \quad \xi^t=\Delta \xi^t=0 \quad \text { on } \quad \partial \Omega \times \mathbb{R}^{+},\\
			w(x, 0)=w_0(x), \quad w_t(x, 0)=w_1(x), \quad \xi^t(x, 0)=0, \quad \xi^0(x, s)=\xi_0(x, s).
		\end{array}\right.
	\end{equation}

	Set
	\begin{equation}\label{eq4.4}
		E_w(t)=\left\|w_t\right\|^2+\left\|\nabla w_t\right\|^2+\|\Delta w\|^2+\left\|\xi^t\right\|_{\mu, 2}^2.
	\end{equation}

Due to the growth condition of the nonlinearity $f(u)$, we will discuss the following in two cases.

\textbf{Case I:}\quad $\left|f^{\prime}(s)\right| \leq M_f\left(1+|s|^p\right),\quad 0<p<p^*$.

 At first, by testing equation $(\ref{eq4.3})_1$ with $w_t$ and $(\ref{eq4.3})_2$ with $\xi^t$, followed by integration over the spatial domain $\Omega$, one obtains
 \begin{equation}\label{eq4.5}
 	\begin{aligned}
 		 \frac{1}{2} \frac{\mathrm{~d}}{\mathrm{~d} t} E_w(t)=&-\|\nabla w_t\|^2-\left(\xi^t, \xi_s^t\right)_{\mu, 2}-(a_{\epsilon}(t)w_t,w_t)\\&-\left(f(u_1)-f(u_2), w_t(t)\right)+(g_1(t)-g_2(t),w_t(t)).
 	\end{aligned}
 \end{equation}

 We now estimate the right-hand side of the (\ref{eq4.5}).
 Similar to Theorem \ref{th3.1}, we have
 \begin{equation}\label{eq4.6}
 	-\left(\xi^t, \xi_s^t\right)_{\mu, 2} \leq -\frac{\delta}{2}\left\|\xi^t\right\|_{\mu, 2}^2,
 \end{equation}
 \begin{equation}\label{eq4.7}
 	-(a_{\epsilon}(t)w_t,w_t)\leq -a_0\|w_t\|^2.
 \end{equation}
 
 Further more, since $\frac{p}{2(p+1)}+\frac{1}{2(p+1)}+\frac{1}{2}=1$, again by generalized H$\ddot{o}$lder inequality and (\ref{eq2.3}) it follows that
 \begin{equation}\label{eq4.8}
 	\begin{aligned}
 		& \left|\left(f(u_1)-f(u_2), w_t\right)\right| \\
 		\leq & M_f \int_{\Omega}\left(1+|u_1|^p+|u_2|^p\right)|w|\left|w_t\right| \mathrm{d} x \\
 		\leq & M_f\left(|\Omega|^{\frac{p}{2(p+1)}}+\|u_1\|_{2(p+1)}^p+\|u_2\|_{2(p+1)}^p\right)\|w\|_{2(p+1)}\left\|w_t\right\| \\
 		\leq & C_B\|w\|_{2(p+1)}\left\|\nabla w_t\right\|\\
 		\leq& \frac{\lambda_0C_B^2}{4a_0}\|w\|_{2(p+1)}^2+a_0\|w_t\|^2.
 	\end{aligned}
 \end{equation}
 Letting us denote $\bar{C_B}=\frac{\lambda_0C_B^2}{4a_0}$, then
 \begin{equation}\label{eq4.9}
 	\left|\left(f(u_1)-f(u_2), w_t\right)\right|\leq \bar{C_B}\|w\|_{2(p+1)}^2+a_0\|w_t\|^2.
 \end{equation}
 
 Inserting (\ref{eq4.6})--(\ref{eq4.9}) into (\ref{eq4.5}), we obtain
 \begin{equation}\label{eq4.10}
 	\frac{d}{dt}E_w(t)\leq -2\|\nabla w_t\|^2-\delta\|\xi^t\|^2_{\mu,2}+\bar{C_B}\|w\|_{2(p+1)}+(g_1(t)-g_2(t),w_t).
 \end{equation}
 
 Next, define the perturbation of (\ref{eq4.4})
 \begin{equation}\label{eq4.11}
 	E^{\varepsilon}_w(t)=E_w(t)+\varepsilon \Phi(t),
 \end{equation}
 	where
 \begin{equation}\label{eq4.12}
 	\Phi(t)=\int_{\Omega} w_t(t) w(t) \mathrm{d} x-\int_{\Omega} \triangle w_t(t) w(t) \mathrm{d} x.
 \end{equation}
 	
  There exists a constant $C_5>0$ directly resulting from the definitions of $E_w(t)$ and $\Phi(t)$ such that
 \begin{equation}\label{eq4.13}
 	\left|E_w^{\epsilon}(t)-E_w(t)\right| \leq \epsilon C_5 E_w(t),
 \end{equation}
 and
 \begin{equation}\label{eq4.14}
 		\begin{aligned}
 		\frac{d}{dt}\Phi(t)&=\int_{\Omega}\left(w_{t t}-\Delta w_{t t}\right) w \mathrm{~d} x+\left\|w_t\right\|_2^2+\left\|\nabla w_t\right\|_2^2\\
 		& =-\frac{1}{2}E_w(t)+\frac{3}{2}\|w_t\|^2+\frac{3}{2}\|\nabla w_t\|^2-\frac{1}{2}\|\Delta w\|^2+\frac{1}{2}\left\|\xi^t\right\|_{\mu, 2}^2-\int_0^{\infty} \mu(s)\left(\Delta \xi^t(s), \Delta w(t)\right) \mathrm{d} s\\
 		&\quad-\int_{\Omega}a_{\epsilon}(t)w_twdx -\int_{\Omega} \nabla w_t \cdot \nabla w \mathrm{~d} x-(f(u_1)-f(u_2),w)+(g_1(t)-g_2(t),w) \mathrm{d} x.
 	\end{aligned}
 \end{equation}

 Let us estimate the right-hand side of (\ref{eq4.14}).
 \begin{equation}\label{eq4.15}
 	\left|\int_0^{\infty} \mu(s)\left(\Delta \xi^t(s), \Delta w(t)\right) \mathrm{d} s\right|\leq \epsilon_2 \mu_0\|\Delta w(t)\|^2+\frac{1}{4 \epsilon_2}\left\|\xi^t\right\|_{\mu, 2}^2,
 \end{equation}
 \begin{equation}\label{eq4.16}
 	\left|\int_{\Omega}a_{\epsilon}(t)w_twdx\right|\leq \frac{\epsilon_2}{\lambda_2}\|\Delta w(t)\|_2^2+\frac{1}{4 \epsilon_2}\left\|\nabla w_t(t)\right\|_2^2,
 \end{equation}
 \begin{equation}\label{eq4.17}
 	\left|\int_{\Omega}a_{\epsilon}(t)w_twdx\right|\leq \frac{a_1^2\epsilon_2}{\lambda_1}\|\Delta w\|^2+\frac{1}{4\lambda_0\epsilon_2}\|\nabla w_t\|^2,
 \end{equation}
 \begin{equation}\label{eq4.18}
 	\left|(f(u_1)-f(u_2),w)\right|\leq C_B\|w\|^2_{2(p+1)}.
 \end{equation}
 
 Inserting (\ref{eq4.15}--(\ref{eq4.18}) into (\ref{eq4.14}), we arrive at
 	\begin{equation}\label{eq4.19}
 		\begin{aligned}
 			\frac{d}{dt}\Phi(t) \leq & -\frac{1}{2} E_w(t)+\left(\frac{3}{2 \lambda_0}+\frac{3}{2}+\frac{1}{4 \epsilon_2}+\frac{1}{4\lambda_0\epsilon_2}\right)\left\|\nabla w_t(t)\right\|_2^2+\left(\frac{1}{2}+\frac{1}{4\epsilon_2}\right)\|\xi^t\|^2_{\mu,2}\\
 			&+\left[\frac{1}{2}-\epsilon_2\left(\mu_0+\frac{a_1^2}{\lambda_1}+\frac{1}{\lambda_2}\right)\right]\|\Delta w\|^2+C_B\|w\|^2_{2(p+1)}+(g_1(t)-g_2(t),w).
 		\end{aligned}
 	\end{equation} 	

Denoting $C_6=\frac{3}{2 \lambda_0}+\frac{3}{2}+\frac{1}{4 \epsilon_2}+\frac{1}{4\lambda_0\epsilon_2}$, $C_7=\frac{1}{2}+\frac{1}{4\epsilon_2}$ and choosing $\epsilon_2$ small enough, such that $\frac{1}{2}-\epsilon_2\left(\mu_0+\frac{a_1^2}{\lambda_1}+\frac{1}{\lambda_2}\right)>0$, then we conclude
\begin{equation}\label{eq4.20}
	\frac{d}{dt}\Phi(t) \leq-\frac{1}{2}E_w(t)+C_6\|\nabla w_t\|^2+C_7\|\xi^t\|^2_{\mu,2}+C_B\|w\|^2_{2(p+1)}+(g_1(t)-g_2(t),w).
\end{equation}
 
From (\ref{eq4.10}) and (\ref{eq4.20}), there exists $\varepsilon_3=\min\{\frac{2}{C_6},\frac{\delta}{C_7},1\}$, then for any $\bar{C_B}$, we deduce
\begin{equation}\label{eq4.21}
	\frac{d}{dt}E_w^{\varepsilon}(t)\leq -\frac{\varepsilon}{2}E_w(t)+2\bar{C_B}\|\nabla w\|^2_{2(p+1)}+(g_1(t)-g_2(t),w_t)+(g_1(t)-g_2(t),w),
\end{equation}
for any $t\geq0$, $\varepsilon\in(0,\varepsilon_3]$.

 Now we take $\varepsilon_4=\min \left\{\frac{1}{2 C_5},\varepsilon_3\right\}$, and choose $\varepsilon \leq \varepsilon_4$. Then (\ref{eq4.13}) implies that
 \begin{equation}\label{eq4.22}
 	 \frac{1}{2} E_w(t) \leq E_w^{\varepsilon}(t) \leq \frac{3}{2} E_w(t), \quad t \geq 0.
 \end{equation}

From (\ref{eq4.21}) and (\ref{eq4.22}), we obtain
\begin{equation}\label{eq4.23}
	\begin{aligned}
		E_w^{\varepsilon}(T)&\leq E_w^{\varepsilon}(s)e^{-\frac{\varepsilon}{3}(T-s)}+2\bar{C_B}\int_{s}^{T}e^{-\frac{\varepsilon}{3}(T-t)}\|w(t)\|^2_{2(p+1)}dt\\
		&\quad+\int_{s}^{T}\int_{\Omega}e^{-\frac{\varepsilon}{3}(T-t)}(g_1(t)-g_2(t))w(t)dxdt\\
		&\quad+\int_{s}^{T}\int_{\Omega}e^{-\frac{\varepsilon}{3}(T-t)}(g_1(t)-g_2(t))w_t(t)dxdt, \quad 0\leq s\leq T.
	\end{aligned}
\end{equation}
 
 Combining (\ref{eq4.22}) and (\ref{eq4.23}), we have
 \begin{equation}\label{eq4.24}
 	\begin{aligned}
 		E_w(T)&\leq 3E_w(s)e^{-\frac{\varepsilon}{3}(T-s)}+4\bar{C_B}\int_{s}^{T}e^{-\frac{\varepsilon}{3}(T-t)}\|w(t)\|^2_{2(p+1)}dt\\
 		&\quad+2\int_{s}^{T}\int_{\Omega}e^{-\frac{\varepsilon}{3}(T-t)}(g_1(t)-g_2(t))w(t)dxdt\\
 		&\quad+2\int_{s}^{T}\int_{\Omega}e^{-\frac{\varepsilon}{3}(T-t)}(g_1(t)-g_2(t))w_t(t)dxdt, \quad 0\leq s\leq T.
 	\end{aligned}
 \end{equation}

 Thanks to \cite{YX15}, similar to the proof of Theorem \ref{th3.1}, we can directly get $E_w(t)\leq \rho_1$, where $\rho_1$ is a positive constant depending on $\|(u_0^i,u_1^i,\eta^{0,i})\|_{\mathcal{H}},i=1,2$, then
 $$3E_w(s)e^{-\frac{\varepsilon}{3}(T-s)}\leq 3\rho_1e^{-\frac{\varepsilon}{3}(T-s)}.$$
  For given $\varepsilon_T>0$, we can choose $T=T(\varepsilon_T,B_0)>s+1$, such that $3\rho_1e^{-\frac{\varepsilon}{3}(T-s)}\leq \varepsilon_T$, to deduce that
 \begin{equation}\label{eq4.25}
 	\begin{aligned}
 		E_w(T)&\leq \varepsilon_T+4\bar{C_B}\int_{s}^{T}e^{-\frac{\varepsilon}{3}(T-t)}\|w(t)\|^2_{2(p+1)}dt\\
 		&\quad+2\int_{s}^{T}\int_{\Omega}e^{-\frac{\varepsilon}{3}(T-t)}(g_1(t)-g_2(t))w(t)dxdt\\
 		&\quad+2\int_{s}^{T}\int_{\Omega}e^{-\frac{\varepsilon}{3}(T-t)}(g_1(t)-g_2(t))w_t(t)dxdt, \quad 0\leq s\leq T.
 	\end{aligned}
 \end{equation}
 
 Integrating (\ref{eq4.25}) over $[0,T]$ with respect to $s$, we obtain that
 \begin{equation}\label{eq4.26}
 	\begin{aligned}
 		E_w(T)&\leq \varepsilon_T+4\bar{C_B}\int_{0}^{T}e^{-\frac{\varepsilon}{3}(T-t)}\|w(t)\|^2_{2(p+1)}dt\\
 		&\quad+\frac{2}{T}\int_{0}^{T}\int_{s}^{T}\int_{\Omega}e^{-\frac{\varepsilon}{3}(T-t)}(g_1(t)-g_2(t))w(t)dxdtds\\
 		&\quad+\frac{2}{T}\int_{0}^{T}2\int_{s}^{T}\int_{\Omega}e^{-\frac{\varepsilon}{3}(T-t)}(g_1(t)-g_2(t))w_t(t)dxdtds, \quad 0\leq s\leq T.
 	\end{aligned}
 \end{equation}
 
 Set 
 \begin{equation*}
 	\begin{aligned}
 			\phi_T\left(y_1,y_2; \sigma_1, \sigma_2\right)&=4\bar{C_B}\int_{0}^{T}e^{-\frac{\varepsilon}{3}(T-t)}\|w(t)\|^2_{2(p+1)}dt\\
 			&\quad+\frac{2}{T}\int_{0}^{T}\int_{s}^{T}\int_{\Omega}e^{-\frac{\varepsilon}{3}(T-t)}(g_1(t)-g_2(t))w(t)dxdtds\quad\\
 			&\quad+\frac{2}{T}\int_{0}^{T}\int_{s}^{T}\int_{\Omega}e^{-\frac{\varepsilon}{3}(T-t)}(g_1(t)-g_2(t))w_t(t)dxdtds, \quad 0\leq s\leq T.
 	\end{aligned}
 \end{equation*}
 Finally, we obtain (\ref{eq4.2}).

\textbf{Case II:}\quad $\left|f^{\prime}(s)\right| \leq M_f\left(1+|s|^p\right),\quad p=p^*$.

At first, multiplying $(\ref{eq4.3})_1$ by $w_t$ and $(\ref{eq4.3})_2$ by $\xi^t$, integrating over $[r,T]
\times\Omega,\,0\leq r\leq T$ and using (\ref{eq4.6}), (\ref{eq4.7}), we deduce that
\begin{equation}\label{eq4.27}
	\begin{aligned}
		&\frac{1}{2}E_w(T)+\int_{r}^{T}\int_{\Omega}|\nabla  w_t(\tau)|^2dxd\tau+\frac{\delta}{2}\int_{r}^{T}\|\xi^t\|^2_{\mu,2}d\tau+a_0\int_{r}^{T}\int_{\Omega}|w_t(\tau)|^2dxd\tau\\
		\leq& \frac{1}{2}E_w(r)+\int_{r}^{T}\int_{\Omega}(g_1(\tau)-g_2(\tau))w_t(\tau)dxd\tau-\int_{r}^{T}\int_{\Omega}(f(u_1(\tau)-f(u_2(\tau)))w_t(\tau)dxd\tau.
	\end{aligned}
\end{equation}

Choosing $C_{\delta,a_0}=\max\{1,\frac{\delta}{2},a_0\}$, then
\begin{equation}\label{eq4.28}
	\begin{aligned}
		&\int_{r}^{T}\int_{\Omega}|\nabla  w_t(\tau)|dxd\tau+\int_{r}^{T}\|\xi^t\|^2_{\mu,2}d\tau+a_0\int_{r}^{T}\int_{\Omega}|w_t(\tau)|^2dxd\tau\\
		\leq& \frac{1}{2}C_{\delta,a_0}E_w(r)+C_{\delta,a_0}\int_{r}^{T}\int_{\Omega}(g_1(\tau)-g_2(\tau))w_t(\tau)dxd\tau\\
		&-C_{\delta,a_0}\int_{r}^{T}\int_{\Omega}(f(u_1(\tau)-f(u_2(\tau)))w_t(\tau)dxd\tau.
	\end{aligned}
\end{equation}

Multiplying (\ref{eq4.3}) by $w_t$ and integrating over $[0,T]\times \Omega$, we deduce
\begin{equation}\label{eq4.29}
	\begin{aligned}
		&\int_{0}^{T}\int_{\Omega}|\Delta w(r)|^2dxdr\\
		\leq&-\int_{\Omega}w(T)w_t(T)dx+\int_{\Omega}w(0)w_t(0)dx+\int_{0}^{T}\int_{\Omega}|w_t(r)|^2dxdr-\int_{0}^{T}\|\xi^t\|^2_{\mu,2}dr\\
		&-\int_{\Omega}\nabla w(T)\nabla w_t(T)dx+\int_{\Omega}\nabla w(0)\nabla w_t(0)dx-\frac{a_0}{2}\int_{\Omega}|w(T)|^2dx+\frac{a_0}{2}\int_{\Omega}|w(0)|^2dx\\
		&+\int_{0}^{T}\int_{\Omega}(g_1(r)-g_2(r))w(r)dxdr-\int_{0}^{T}\int_{\Omega}(f(u_1(r))-f(u_2(r)))w(r)dxdr.
	\end{aligned}
\end{equation}
From (\ref{eq4.28}) and (\ref{eq4.29}), we derive
\begin{equation}\label{eq4.30}
	\begin{aligned}
		&\int_{0}^{T}E_w(r)dr
		\leq-\int_{\Omega}w(T)w_t(T)dx+\int_{\Omega}w(0)w_t(0)dx+\int_{0}^{T}\int_{\Omega}|w_t(r)|^2dxdr-\int_{0}^{T}\|\xi^t\|^2_{\mu,2}dr\\
		&\quad-\int_{\Omega}\nabla w(T)\nabla w_t(T)dx+\int_{\Omega}\nabla w(0)\nabla w_t(0)dx-\frac{a_0}{2}\int_{\Omega}|w(T)|^2dx+\frac{a_0}{2}\int_{\Omega}|w(0)|^2dx\\
		&\quad+\int_{0}^{T}\int_{\Omega}(g_1(r)-g_2(r))w(r)dxdr-\int_{0}^{T}\int_{\Omega}(f(u_1(r))-f(u_2(r)))w(r)dxdr+C_{\delta,a_0}E_w(r)\\
		&\quad+2C_{\delta,a_0}\int_{r}^{T}\int_{\Omega}(g_1(\tau)-g_2(\tau))w_t(\tau)dxd\tau-2C_{\delta,a_0}\int_{r}^{T}\int_{\Omega}(f(u_1(\tau)-f(u_2(\tau)))w_t(\tau)dxd\tau.
	\end{aligned}
\end{equation}

Integrating (\ref{eq4.27}) over $[0,T]$ with respect to $r$ and using (\ref{eq4.30}), we obtain that
\begin{equation}\label{eq4.31}
	\begin{aligned}
		&TE_w(T)\\
		\leq&2\int_{0}^{T}\int_{r}^{T}\int_{\Omega}(g_1(\tau)-g_2(\tau))w_t(\tau)dxd\tau\\
		&-2\int_{0}^{T}\int_{r}^{T}\int_{\Omega}(f(u_1(\tau)-f(u_2(\tau)))w_t(\tau)dxd\tau+\int_{0}^{T}E_w(r)dr\\
		\leq&2\int_{0}^{T}\int_{r}^{T}\int_{\Omega}(g_1(\tau)-g_2(\tau))w_t(\tau)dxd\tau-2\int_{0}^{T}\int_{r}^{T}\int_{\Omega}(f(u_1(\tau)-f(u_2(\tau)))w_t(\tau)dxd\tau\\
		&-\int_{\Omega}w(T)w_t(T)dx+\int_{\Omega}w(0)w_t(0)dx-\int_{\Omega}\nabla w(T)\nabla w_t(T)dx+\int_{\Omega}\nabla w(0)\nabla w_t(0)dx\\
		&-\frac{a_0}{2}\int_{\Omega}|w(T)|^2dx+\frac{a_0}{2}\int_{\Omega}|w(0)|^2dx+\int_{0}^{T}\int_{\Omega}(g_1(r)-g_2(r))w(r)dxdr\\
		&-\int_{0}^{T}\int_{\Omega}(f(u_1(r))-f(u_2(r)))w(r)dxdr+C_{\delta,a_0}E_w(0)\\
		&+2C_{\delta,a_0}\int_{r}^{T}\int_{\Omega}(g_1(\tau)-g_2(\tau))w_t(\tau)dxd\tau\\
		&-2C_{\delta,a_0}\int_{r}^{T}\int_{\Omega}(f(u_1(\tau)-f(u_2(\tau)))w_t(\tau)dxd\tau.
	\end{aligned}
\end{equation}

Set
 \begin{equation*}
	\begin{aligned}
		\phi_T\left(y_1,y_2; \sigma_1, \sigma_2\right)&=2\int_{0}^{T}\int_{r}^{T}\int_{\Omega}(g_1(\tau)-g_2(\tau))w_t(\tau)dxd\tau\\
		&\quad-2\int_{0}^{T}\int_{r}^{T}\int_{\Omega}(f(u_1(\tau)-f(u_2(\tau)))w_t(\tau)dxd\tau\\
		& \quad+\int_{0}^{T}\int_{\Omega}(g_1(r)-g_2(r))w(r)dxdr-\int_{0}^{T}\int_{\Omega}(f(u_1(r))-f(u_2(r)))w(r)dxdr\\
		&\quad+2C_{\delta,a_0}\int_{r}^{T}\int_{\Omega}(g_1(\tau)-g_2(\tau))w_t(\tau)dxd\tau\\
		&\quad-2C_{\delta,a_0}\int_{r}^{T}\int_{\Omega}(f(u_1(\tau)-f(u_2(\tau)))w_t(\tau)dxd\tau\quad 0\leq s\leq T,
	\end{aligned}
\end{equation*}
\begin{equation*}
	\begin{aligned}
		C_M&=-\int_{\Omega}w(T)w_t(T)dx+\int_{\Omega}w(0)w_t(0)dx-\int_{\Omega}\nabla w(T)\nabla w_t(T)dx+\int_{\Omega}\nabla w(0)\nabla w_t(0)dx\\
		&\quad-\frac{a_0}{2}\int_{\Omega}|w(T)|^2dx+\frac{a_0}{2}\int_{\Omega}|w(0)|^2dx.
	\end{aligned}
\end{equation*}

Given that the family of processes $\{U_\sigma(t, \tau)\}_{\sigma \in \Sigma}$ admits a bounded absorbing set uniformly in $\sigma$, and invoking Theorem~\ref{th3.1}, it follows that for any prescribed $\varepsilon_T > 0$, there exists a sufficiently large time $T > 0$ such that
$$
\frac{C_M}{T} \leqslant \varepsilon_T,
$$
then we obtain (\ref{eq4.2}). The proof is hence complete.
\end{Lemma}
$\hfill\qedsymbol$

\textbf{\large {4.2 Uniformly asymptotic compactness}}

By virtue of Lemma \ref{le4.1}, it suffices to verify that $\phi_T(\cdot, \cdot ; \cdot, \cdot)$ belongs to the class $\operatorname{Contr}(B_0, \Sigma)$ for each fixed $T > 0$. Following the argument developed in~\cite{CS6}, one readily obtains the following conclusion
\begin{Lemma}\rm(\cite{CS6})\label{le4.2}
	If $g \in L^2\left(\mathbb{R} ; L^2\right)$, and $g_t \in L_b^r\left(\mathbb{R} ; L^r\right)$ with $\gamma \geq \frac{6}{5}$, then there exists a constant $M>0$ such that $\sup\limits _{t \in \mathbb{R}}\|g(t+s)\|_2 \leq M, \quad$ for all $\quad s \in \mathbb{R}$.
\end{Lemma}
\begin{Lemma}\rm(\cite{CS6})\label{le4.3}
	Let $\left\{s_n\right\}_{n \in \mathbb{N}} \subset \mathbb{R}$ and $\left\{u_n\right\}_{n \in \mathbb{N}}$ be a sequence uniformly bounded in $L^{\infty}\left(\mathbb{R}^{+} ; H_0^1(\Omega)\right)$ such that: for any $T_1>0,\left\{u_{n t}\right\}_{n \in \mathbb{N}}$ is uniformly bounded in $L^{\infty}\left(0, T_1 ; L^2(\Omega)\right)$. If $g$ verifies (\ref{eq2.8})--(\ref{eq2.10}), then for any $T>0$, there exists subsequences $\left\{u_{n_k}\right\}_{k \in \mathbb{N}} \subset\left\{u_n\right\}_{n \in \mathbb{N}}$ and $\left\{s_{s_k}\right\}_{k \in \mathbb{N}} \subset\left\{s_n\right\}_{n \in \mathbb{N}}$ such that
	$$
	\lim _{k \rightarrow \infty} \lim _{l \rightarrow \infty} \int_0^T \int_s^T \int_{\Omega}\left(g\left(x, \tau+s_{n_k}\right)-g\left(x, \tau+s_{n_l}\right)\right)\left(u_{n_k}-u_{n_l}\right)_t(\tau) d x d \tau d s=0.
	$$
\end{Lemma}
\begin{Theorem}\label{th4.4}
Assuming that conditions (\ref{eq2.3})--(\ref{eq2.13}) are satisfied, and that the symbol space $\Sigma$ is given by~(\ref{eq2.9}), it follows that the family of processes $\{U_\sigma(t, \tau)\}_{\sigma \in \Sigma}$ associated with system (\ref{eq1.5})--(\ref{eq1.7}) is asymptotically compact in $\mathcal{H}$ uniformly with respect to $\sigma \in \Sigma$.
 
 \rm\textbf{Proof.} At first, based on the argument established in the proof of Theorem~\ref{th3.1}, it can be concluded that for any fixed $T > 0$,
 \begin{equation}\label{eq4.32}
 	\bigcup_{\sigma \in \Sigma} \bigcup_{t \in[0, T]} U_\sigma(t, 0) B_0 \text { is bounded in } E_0,
 \end{equation}
 and the bound depends on $T$.
 
 \textbf{Case I:}\quad $\left|f^{\prime}(s)\right| \leq M_f\left(1+|s|^p\right),\quad 0<p<p^*$.
 
 We first denote 
 \begin{equation}\label{eq4.33}
 	\Phi_T\left(y^1, y^2 ; \sigma_1, \sigma_2\right)=\Phi_T^1\left(y^1, y^2\right)+\Phi_T^2\left(\sigma_1, \sigma_2\right),
 \end{equation}
where
 $$
 \Phi_T^1\left(y^1, y^2\right)=4\bar{C_B}\int_{0}^{T}e^{-\frac{\varepsilon}{3}(T-t)}\|w(t)\|^2_{2(p+1)}dt,
 $$
 and
 \begin{equation*}
 	\begin{aligned}
 		\Phi_T^2\left(\sigma^1,\sigma^2\right)&=\frac{2}{T}\int_{0}^{T}\int_{s}^{T}\int_{\Omega}e^{-\frac{\varepsilon}{3}(T-t)}(g_1(t)-g_2(t))w(t)dxdtds\\
 		&\quad+\frac{2}{T}\int_{0}^{T}\int_{s}^{T}\int_{\Omega}e^{-\frac{\varepsilon}{3}(T-t)}(g_1(t)-g_2(t))w_t(t)dxdtds, \quad 0\leq s\leq T.
 	\end{aligned}
 \end{equation*}
 Since $H^2(\Omega) \cap H_0^1(\Omega) \hookrightarrow L^{2(p+1)}(\Omega)$ and from \cite{YX16}, we infer that there exists a subsequence $\left\{u_{n_k}\right\}_{k \in \mathbb{N}}$ which converges strongly in $C([0, T]$; $\left.L^2(\Omega)\right)$ such that
 \begin{equation}\label{eq4.34}
 \liminf _{k \rightarrow \infty} \liminf _{l \rightarrow \infty} \Phi_T^1\left(\left(u_{n_k}^0, u_{n_k}^1\right) ;\left(u_{n_l}^0, u_{n_l}^1\right)\right)=0 .
 \end{equation}
 Thanks to Lemma \ref{le4.2} and Lemma \ref{le4.3}, we can directly get
 	\begin{equation}\label{eq4.35}
 		\liminf _{k \rightarrow \infty} \liminf _{l \rightarrow \infty} \Phi_T^2\left(\sigma_{n_k} ; \sigma_{n_l}\right)=0
 	\end{equation}
Hence, combining (\ref{eq4.33})-(\ref{eq4.35}), we get $\phi_T(\cdot, \cdot ; \cdot, \cdot) \in \operatorname{Contr}\left(B_0, \Sigma\right)$ immediately.

\textbf{Case II:}\quad $\left|f^{\prime}(s)\right| \leq M_f\left(1+|s|^p\right),\quad p=p^*$.

Again by using $H^2(\Omega) \cap H_0^1(\Omega) \hookrightarrow L^{2(p+1)}(\Omega)$ and Lemma \ref{le4.2}, similar to the method used in Lemma \ref{le4.3}, we get
\begin{equation}\label{eq4.36}
	\lim _{n \rightarrow \infty} \lim _{m \rightarrow \infty} \int_0^T \int_{\Omega}\left(g_n(x, r)-g_m(x, r)\right)\left(u_{n}(r)-u_{m}(r)\right) d x d r=0,
\end{equation}
\begin{equation}\label{eq4.37}
	\lim _{n \rightarrow \infty} \lim _{m \rightarrow \infty} \int_0^T \int_{\Omega}\left(g_n(x, r)-g_m(x, r)\right)\left(u_{n_t}(r)-u_{m_t}(r)\right) d x d r=0,
\end{equation}
\begin{equation}\label{eq4.38}
	\lim _{n \rightarrow \infty} \lim _{m \rightarrow \infty} \int_0^T \int_r^T \int_{\Omega}\left(g_n(x, \tau)-g_m(x, \tau)\right)\left(u_{n_t}(\tau)-u_{m_t}(\tau)\right) d x d \tau d r=0.
\end{equation}
At the same time, from \cite{CS6} and the growth condition (\ref{eq2.3})--(\ref{eq2.5}), we can easily get that
\begin{equation}\label{eq4.39}
	\lim _{n \rightarrow \infty} \lim _{m \rightarrow \infty} \int_0^T \int_{\Omega}\left(f\left(u_n(r)\right)-f\left(u_m(r)\right)\right)\left(u_n(r)-u_m(r)\right) d x d r=0,
\end{equation}
\begin{equation}\label{eq4.40}
	\lim _{n \rightarrow \infty} \lim _{m \rightarrow \infty} \int_0^T \int_{\Omega}\left(f\left(u_n(r)\right)-f\left(u_m(r)\right)\right)\left(u_{n_t}(r)-u_{m_t}(r)\right) d x d r=0,
\end{equation}
\begin{equation}\label{eq4.41}
	\lim _{n \rightarrow \infty} \lim _{m \rightarrow \infty} \int_0^T\int_{r}^{T} \int_{\Omega}\left(f\left(u_n(\tau)\right)-f\left(u_m(\tau)\right)\right)\left(u_{n_t}(\tau)-u_{m_t}(\tau)\right) d xd\tau d r=0.
\end{equation}

Hence, combining (\ref{eq4.36})-(\ref{eq4.41}), we get that $\phi_{T}(\cdot, \cdot ; \cdot, \cdot) \in \operatorname{Contr}\left(B_0, \Sigma\right)$, and then this completes the proof of Theorem \ref{th4.4}.
\end{Theorem}
 $\hfill\qedsymbol$

\begin{Theorem}\label{th4.5}
Suppose that assumptions (\ref{eq2.3})--(\ref{eq2.13}) hold, and let $\Sigma$ be the symbol space defined in~(\ref{eq2.9}). Then, the family of processes $\{U_\sigma(t, \tau)\}_{\sigma \in \Sigma}$ associated with problem~(\ref{eq1.5})--(\ref{eq1.7}) admits a compact uniform attractor $\mathcal{A}_{\Sigma}$ in $\mathcal{H}$, which is uniform with respect to $\sigma \in \Sigma$.
	
\rm\textbf{Proof.} From Theorem \ref{th2.9} , we know Theorems \ref{th3.1} and \ref{th4.4} imply the existence of an uniform attractor immediately.
\end{Theorem}
$\hfill\qedsymbol$

\section{Upper semicontinuity of the uniform attractors}

 \quad In this section, we aim to demonstrate the upper semicontinuity of uniform attractors relative to the functional parameter $a$. To this end, let $\{a_{\epsilon} : \epsilon \in [0,1]\}$ be a family of real-valued functions satisfying condition (\ref{eq2.6}). Denote by $\{U_{a_{\epsilon}, \sigma}(\cdot, \cdot)\}_{\sigma \in \Sigma}$ the corresponding evolution processes, and by $\{\mathcal{A}_{\Sigma}^{a_{\epsilon}}\}_{\epsilon \in [0,1]}$ the associated family of uniform attractors for problem~(\ref{eq1.5})--(\ref{eq1.7}).
  
  Additionally, suppose that
  \begin{equation}\label{eq5.1}
  	\left\|a_{\epsilon}-a_0\right\|_{L^{\infty}(\mathbb{R})} \rightarrow 0 \text { as } \epsilon \rightarrow 0^{+}.
  \end{equation}
  
  Now we present the main theorem that guarantee the uppersemicontinuity of uniform attractors as follows.
  \begin{Theorem}\label{th5.1}
   Assume that the conditions stated in Theorem \ref{th4.5} are satisfied. For each $\epsilon \in [0,1]$, let $w^{(\epsilon)}(\cdot) := U^{a_{\epsilon}}_{\Sigma}(\cdot, \tau) w_0$ denote the solution to problem~(\ref{eq1.5})--(\ref{eq1.7}) in the phase space $\mathcal{H}$. Then, for every $T > 0$, one has
   \[
   w^{(\epsilon)} \to w^{(0)} \quad \text{in } C([0, T]; \mathcal{H}) \quad \text{as } \epsilon \to 0^+.
   \]
   Furthermore, the family of uniform attractors $\{\mathcal{A}^{a_{\epsilon}}_{\Sigma}\}_{\epsilon \in [0,1]}$ exhibits upper semicontinuity at $\epsilon = 0$., that is
  	\begin{equation}\label{eq5.2}
  			\lim _{\epsilon \rightarrow 0} \operatorname{dist}_{\mathcal{H}}\left(A^{a_\epsilon}_{\Sigma}, A^{a_0}_{\Sigma}\right)=0,
  	\end{equation}
  and
  \begin{equation}\label{eq5.2-1}
  	\lim _{\epsilon \rightarrow 0} \operatorname{dist}_{\mathcal{H}}\left(\mathcal{K}_g^{a_\epsilon}(s), \mathcal{K}_g^0(s)\right)=0, \quad \forall s \in \mathbb{R}, g \in \Sigma,
  \end{equation}
where $\mathcal{K}_g^{a_\epsilon}$ is the kernel of the process $U_g^{a_\epsilon}, \mathcal{K}_g^{a_\epsilon}(s)$ is the kernel section at time $s, \epsilon\in [0,1]$, and $dist_{\mathcal{H}}$ denotes the standard Hausdorff semidistance $\operatorname{in} \mathcal{H}$.

  \rm\textbf{Proof.} For each $w_0 \in \mathcal{H}$, consider $u=w^{(\epsilon)}=U^{a_{\epsilon} }_{\Sigma}(t, \tau) w_0$ and $v=w^{(0)}=U^{a_0 }_{\Sigma}(t, \tau) w_0$. Let $w=u-v,\,\xi^t=\eta^t-\nu^t$. Then
  \begin{equation}\label{eq5.3}
  	\left\{\begin{array}{l}
  		w_{t t}-\Delta w_{tt}+\Delta^2 w+\int_{0}^{\infty}\mu(s)\Delta^2\xi^t(s)ds=\Delta w_t+f(v)-f(u)+a_\epsilon(t) u_t -a_0(t)v_t,\\
  		\xi^t_t(x,s)=-\xi^t_s(x,s)+u_t(x,s).
  	\end{array}\right.
  \end{equation}
  
  Multiplying $(\ref{eq5.3})_1$ by $w_t$ and $(\ref{eq5.3})_2$ by $\xi_t$, we get
  \begin{equation}\label{eq5.4}
  	\begin{aligned}
  		&\frac{1}{2}\frac{d}{dt}(\|\Delta w\|^2+\|\nabla w_t\|^2+\|w_t\|^2+\|\xi^t\|_{\mu,2}^2)\\
  		=&-(\xi^t,\xi^t_s)_{\mu,2}+\int_{\Omega}\Delta w_tw_tdx+\int_{\Omega}(f(v)-f(u))w_tdx+\int(a_\epsilon(t) u_t -a_0(t)v_t)w_tdx.
  	\end{aligned}
  \end{equation}
  
  Note that
  \begin{equation}\label{eq5.5}
  	\begin{aligned}
  		\int_{\Omega}\left(a_\epsilon(t) u_t-a_0(t) v_t\right) w_t d x & =\int_{\Omega}\left(a_\epsilon(t)-a_0(t)\right) u_t w_t d x+\int_{\Omega} a_0(t) w_t w_t d x \\
  		& \leq\left\|a_\epsilon-a_0\right\|_{L^{\infty}(\mathbb{R})}\left\|u_t\right\|\left\|w_t\right\|+a_0\left\|w_t\right\|^2.
  	\end{aligned}
  \end{equation}
  
  From (\ref{eq4.6})--(\ref{eq4.8}), (\ref{eq5.4}), (\ref{eq5.5}) and condition $H^2(\Omega) \cap H_0^1(\Omega) \stackrel{c}{\hookrightarrow} L^{2(p+1)}(\Omega)$, we know there exist positive constants $K_1$ and $K_B$ large enough such that
  \begin{equation}\label{eq5.6}
  	\begin{aligned}
  		&\frac{d}{dt}(\|\Delta w\|^2+\|\nabla w_t\|^2+\|w_t\|^2+\|\xi^t\|_{\mu,2}^2)\\
  		\leq&K_B(\|\Delta w\|^2+\|\nabla w_t\|^2+\|w_t\|^2+\|\xi^t\|_{\mu,2}^2)+K_1\left\|a_\epsilon-a_0\right\|_{L^{\infty}(\mathbb{R})},
  	\end{aligned}
  \end{equation}
  and consequently there exists a positive constant $K_2$
 such that
 \begin{equation}\label{eq5.7}
 	\|\Delta w\|^2+\|\nabla w_t\|^2+\|w_t\|^2+\|\xi^t\|_{\mu,2}^2\leq K_2\left\|a_\epsilon-a_0\right\|_{L^{\infty}(\mathbb{R})}e^{K_Bt},\quad t>0,
 \end{equation}  
  in other words, the solution $u = w^{(\varepsilon)} \left(= U^{a_{\epsilon}}_{\Sigma}(t, \tau) w_0 \right)$ converges to $v = w^{(0)} \left(= U^{a_0}_{\Sigma}(t, \tau) w_0 \right)$ as $\varepsilon \to 0^+$, uniformly on compact subsets of $\mathbb{R}$, with respect to initial data $w_0$ taken from bounded subsets of $\mathcal{H}$.
  
  Now suppose that the family $\mathcal{A}^{a_{\epsilon}}_{\Sigma}$ is not semicontinuous at $\epsilon=0$. Then there exists $\delta^{\prime}>0$ and a sequence $\left\{\epsilon_n\right\} \subset[0, \infty)$ with $\epsilon_n \xrightarrow{n \rightarrow \infty} 0$ such that
  $$
  \operatorname{dist}_{\mathcal{H}}\left(A_{a_{\epsilon_n}, \Sigma}, A_{a_0, \Sigma}\right) \geq \delta^{\prime}.
  $$
  
  Hence, there exists $\left\{y_n\right\}_{n \in \mathbb{N}} \subset \mathcal{A}_{\Sigma}^{a_{\epsilon_n}}$ such that
  \begin{equation}\label{eq5.8}
  	\operatorname{dist}_{\mathcal{H}}\left(y_n, \mathcal{A}_{\Sigma}^{a_{\epsilon_0}}\right) \geq \delta^{\prime}, \quad \forall n \in \mathbb{N} .
  \end{equation}
 
  If $B_0$ is the uniformly (w.r.t. $\sigma \in \Sigma, \epsilon \in[0, \infty)$ ) absorbing set as provided by Theorem \ref{th3.1}. We can select a sufficiently large $m>0$ to ensure that
 \begin{equation}\label{eq5.9}
 	\bigcup_{\epsilon\in[0,1]} \bigcup_{\sigma \in \Sigma} U_\sigma^{a_{\epsilon}}(t, 0) B_0 \subset B_0, \quad \forall t \geq m,
 \end{equation}
 and
 \begin{equation}\label{eq5.10}
 	 \sup _{\sigma\in \Sigma} \operatorname{dist}\left(U_\sigma^{a_{\epsilon_0}}(m, 0) B_0, \mathcal{A}_{\Sigma}^{\epsilon_0}\right) \leq \frac{\delta}{3}^{\prime}.
 \end{equation}

 From Theorem \ref{th2.9} we know $\mathcal{A}_{\Sigma}^{a_{\epsilon_n}}=\omega_{0, \Sigma}^{\epsilon_n}(B_0)(n \in \mathbb{N})$. Therefore, we can identify sequences $\left\{\sigma_n\right\}_{n \in \mathbb{N}} \subset \Sigma,\left\{x_n\right\}_{n \in \mathbb{N}} \subset B_0$ and $\left\{t_n\right\}_{n \in \mathbb{N}} \subset \mathbb{R}^{+}$ with $t_n \rightarrow \infty$, to satisfy, without reducing the generality, that $t_n \geq 2 m$ 
 $$
 \left\|U_{g_n}^{a_{\epsilon_n}}\left(t_n, 0\right) x_n-y_n\right\|_{\mathcal{H}} \leq \frac{\delta^{\prime}}{3}, \quad \forall n \in \mathbb{N}.
 $$
 
 Let $\tilde{x}_n=U_{\sigma_n}^{a_{\epsilon_n}}\left(t_n-m, 0\right) x_n$ and $\sigma_n^{\prime}=T\left(t_n-m\right) \sigma_n$, by (\ref{eq2.16}), (\ref{eq2.17}), (\ref{eq5.9}) and noticing that $t_n \geq$ $2 m$, we have
 $$
 \left\{\tilde{x}_n\right\}_{n \in \mathbb{N}} \subset B_0,
 $$
 and
 $$
 U_{\sigma_n}^{a_{\epsilon_n}}\left(t_n, 0\right) x_n=U_{\sigma_n}^{a_{\epsilon_n}}\left(t_n, t_n-m\right) \widetilde{x}_n=U_{\sigma_n^{\prime}}^{a_{\epsilon_n}}(m, 0) \widetilde{x}_n.
 $$

 On the other hand, due to (\ref{eq5.7}), we can choose $N \in \mathbb{N}$ large enough such that
 $$
 \left\|U_{\sigma_N^{\prime}}^{a_{\epsilon_N}}(m, 0) \widetilde{x}_N-U_{\sigma_N^{\prime}}^{a_{\epsilon_0}}(m, 0) \widetilde{x}_N\right\| \leq \frac{\delta}{3}.
 $$
 
 Therefore, from the above analysis we find
 $$
 \begin{aligned}
 	\operatorname{dist}_{\mathcal{H}}\left(y_N, \mathcal{A}_{\Sigma}^{a_{\epsilon_0}}\right) 
 	\leq& \operatorname{dist}_{\mathcal{H}}\left(y_N, U_{\sigma_N^{\prime}}^{a_{\epsilon_N}}(m, 0) \widetilde{x}_N\right) \\
 	&+\operatorname{dist}_{\mathcal{H}}\left(U_{\sigma_N^{\prime}}^{a_{\epsilon_N}}(m, 0) \widetilde{x}_N, U_{\sigma_N^{\prime}}^{a_{\epsilon_0}}(m, 0) \widetilde{x}_N\right) \\
 	&+\operatorname{dist}_{\mathcal{H}}\left(U_{\sigma_N^{\prime}}^{a_{\epsilon_0}}(m, 0) \widetilde{x}_N, U_{\sigma_N^{\prime}}^{a_{\epsilon_0}}(m, 0) B_0\right) \\
 	&+\operatorname{dist}_{\mathcal{H}}\left(U_{\sigma_N^{\prime}}^{a_{\epsilon_0}}(m, 0) B_0, \mathcal{A}_{\Sigma}^{a_{\epsilon_0}}\right) \\
 	 \leq& \frac{\delta}{3}^{\prime}+\frac{\delta}{3}^{\prime}+0+\frac{\delta}{3}^{\prime}=\delta^{\prime},
 \end{aligned}
 $$
 which contradicts (\ref{eq5.8}). This show that (\ref{eq5.2}) holds and the family of uniform attractors $\mathcal{A}^{a_{\epsilon}}_{\Sigma}$ is upper semicontinuous at $\epsilon=0$.
 
 The proof of (5.3) is analogous to that of (5.2), and the detailed arguments can be found in \cite{YX15}.
  \end{Theorem}
$\hfill\qedsymbol$

$\mathbf{Acknowledgment}$

This work was supported by the TianYuan Special Funds of the NNSF of China with contract number 12226403, the NNSF of China with contract
No.12171082, the fundamental research funds for the central universities with contract numbers CUSF-DH-D-2025032, 2232022G-13, 2232023G-13 and 2232024G-13. 

$\mathbf{Data\,\,availability\,\,statement}$

The data that support the findings of this study are available from the corresponding author upon reasonable request.

$\mathbf{Conflict\,\,of\, \,interest\,\, statement}$ 

No potential conflict of interest was reported by the authors.

\end{document}